%% file: waves_temp_interfaces_arxiv.tex
\documentclass[11pt]{article}
\usepackage{graphics}
\usepackage{graphicx}
\usepackage{pgfplots}
\usepackage{tikz}
\usepackage{amsmath}
\usepackage{mathtools}
\usepackage{epsf}
\usepackage{array}
\usepackage{amsfonts}
\usepackage{amssymb}
\usepackage[english]{babel}
\usepackage[ansinew]{inputenc}              
\usepackage{color}
\usepackage{multirow}
\usepackage{bm}
\usepackage{enumitem}
\usepackage[cal=boondoxo]{mathalfa}
\usepackage{caption}
\usepackage{subcaption}
\usepackage{parskip}
\usepackage{physics}
\usepackage{float}
\begin{document}

\input macro

\title{\bf Edge Waves and Transmissions for Temporal \\ Laminates and Imperfect Chiral Interfaces}

\author{ {\bf A.B. Movchan$^1$, N.V. Movchan$^1$, I.S. Jones$^{1,2}$, }\\ {\bf G.W. Milton$^3$, H.-M. Nguyen$^4$}
\\
\\ $^1$ {\small Department of Mathematical Sciences, University of Liverpool, Liverpool L69 7ZL, U.K.}
\\ $^2$ {\small Liverpool John Moores University, Liverpool L3 3AF, U.K.}
\\ $^3$ {\small  Department of Mathematics, University of Utah, Salt Lake City, Utah 84112, U.S.A.} 
\\ $^4$ {\small Laboratoire Jacques Louis Lions, 
Sorbonne Universit\'{e}, 
4 Place Jussieu, 75013, Paris, France}}

\date{}

\maketitle

\begin{abstract}
The analysis of wave patterns in a structure which possesses periodicity in the spatial and temporal dimensions is presented. 
The topic of imperfect chiral interfaces is also considered. Although causality is fundamental for physical processes, natural wave phenomena can be observed when a wave is split at a temporal interface.   A wave split at a spatial interface is a more common occurrence; however when the coefficients of the governing equations are time-dependent, the temporal interface becomes important.  
Here, the associated edge waves are studied, and  regimes are analysed where the growth of the solution in time is found. 
Imperfect interfaces, across which the displacements are discontinuous, are also considered in the vector case of chiral elastic systems. 
Analytical study and asymptotic approximations are  supplied with illustrative numerical examples. 
\end{abstract}

\section{Introduction}


The phenomenon of wave reflection on spatial interfaces, separating two media, is well-known and is described in classical textbooks (see, for example, books by Kittel \cite{Kittel} and Lekner \cite{Lekner}). If the interfaces are arranged in a finite or infinite array, then multiple reflections may lead to interesting phenomena of dispersion connected to Floquet theory (see \cite{Kittel}). For finite stacks, the book 
\cite{Lekner} has provided constructive algorithms for evaluation of transmission and reflection characteristics of the structured stack.

The emphasis of the present paper is on interfaces of a different kind - temporal interfaces. Interest in this  area has been growing for a number of 
years: see, for example, the book by Lurie \cite{Lurie} and the review of Caloz $\&$ Deck-L\'eger \cite{CD-L}. Most notable has been the work of Fink 
\cite{Fink}, and Bacot {\em et al.} \cite{Fink1} on image recovery through partial time reversal at  temporal interfaces. 
The recent development in modelling of waves in media with temporal and spatial interfaces was stimulated by the work of Lurie \cite{Lurie}, Milton $\&$ Mattei \cite{MilMat}, and Mattei $\&$ Milton \cite{MatMil1, MatMil2}, who introduced the key ideas in this novel subject area.  

With causality in place, a temporal interface still produces wave splitting so that, at a certain time when the elastic or inertial properties of the carrier medium change instantaneously, the wave front will split into two fronts, propagating in different directions. Special combinations of temporal and spatial interfaces can be analysed, where highly non-trivial wave patterns are observed, as discussed in \cite{MilMat, MatMil1, MatMil2}. These patterns of characteristic lines are known as ``field patterns''.

These field patterns can arise in space-time geometries as simple as certain temporal laminates. Here we show that this allows a straightforward computation of their response which, in particular, reveals the behaviour at the edge of the wave. Additionally, the case of chiral interfaces is explored. 
The idea of non-reciprocity for dynamic materials was discussed by  Lurie \cite{Lurie} for the case of a laminate with oblique interfaces in space-time, and by Brun {\em et al.} \cite{Brun2012} in the context of chiral elastic lattices and homogenised coupled elastic systems, This was followed by a computational experiment by Wang {\em et al.} \cite{Bertoldi} who observed a uni-directional edge wave in a homogenised elastic solid. The work on non-reciprocal dynamic materials was further stimulated by novel physical applications and, in particular, by the concept of ``topological insulators'', as discussed by Moore \cite{moore}, Pendry {\em et al.} \cite{pendry}, Hibbins {\em et al.} \cite{hibbins}, and Zhao {\em et al.} \cite{Zhao}.

Floquet theory for transient problems in media with periodic temporal interfaces can be successfully applied, as discussed by Nassar {\em et al.} \cite{Nassar}, and by Lurie $\&$ Weekes \cite{LW}. 

Recent work on modelling of chiral waves in elastic lattices by Carta {\em et al.} \cite{chirallattice2, chirallattice1},   has provided an explanation of the dynamics of discrete gyroscopic systems in the context of dispersion, localisation and dynamic degeneracies.  The paper by Nieves {\em et al.} \cite{GyroMS}
provided mathematical insight into the vibrations of chiral multi-structures and connections between the discrete non-reciprocal systems and their continuous counterparts.  The recent paper by Jones {\em et al.}  \cite{Jones2020} has presented a comprehensive study of coupled elastic waves, dynamic localisation, and dynamic Green's functions in a chiral elastic system. Furthermore, this work laid the foundations for the idea of an imperfect temporal interface associated with high-gyricity material. One characteristic feature of a gyroscopic force is that it is orthogonal to the velocity vector, and hence a rotation transformation is involved in the description of vector chiral waves and also in the analysis of  imperfect temporal chiral interfaces.

The structure of the paper is as follows: In Section \ref{temporal} the formal description of temporal interfaces and temporally stratified media is introduced. For a given set of initial conditions, Section \ref{edgewave} presents the description of the field patterns, 
associated with temporal laminates, and initiates the study of edge waves. 
Special attention is given to the temporal switching algorithm, which leads to the edge wave blow-up. For a different physical configuration related to temporal interfaces, the edge waves were also observed in \cite{LW}.
The case of periodic initial conditions is considered in Section \ref{per_initial}. The notion of chiral vibrations and chiral interfaces is introduced in Section \ref{chiral_interface} and includes the treatment of both spatial and temporal interfaces. Finally, the derivation of the transmission conditions for chiral interfaces, characterised by high gyricity, is presented in Section  \ref{temporal_chiral}. Furthermore, the discussion of the wave splitting and the edge wave is also presented. In that section, the mathematical formulation is set as a vector problem and chiral temporal interfaces provide coupling between the longitudinal and transverse displacements of the wave field.  
Illustrative examples, which include closed form solutions of the Cauchy problems in temporally stratified media, are discussed in Section \ref{S_examples}.

\section{Formal settings - temporal interfaces}
\label{temporal}

Let  two-dimensional vector fields $\Bj$ and  $\Be$  be divergence-free and  curl-free, respectively in $\BR^2$.  Assume that
\beq
\Bj (\Bx)  = \BGs(\Bx) \Be(\Bx),
\eequ{eq1}
where 
\beq
\Grad \cdot \Bj(\Bx) = 0, ~ \Be = - \Grad V,
\eequ{eq2} 
with the diagonal matrix $\BGs (\Bx) = \mbox{diag} (\Ga (\Bx), - \Gb (\Bx)).$ 
Here, $\Ga$ and  $\Gb$ are positive functions. 

Thus, the potential function $V$ satisfies the governing equation
\beq
\fr{\prt}{\prt x_1} \B(\Ga(\Bx) \fr{\prt }{\prt x_1} V(\Bx) \B) - \fr{\prt}{\prt x_2} \B(\Gb(\Bx) \fr{\prt }{\prt x_2} V(\Bx) \B) = 0.
\eequ{eq3}

The initial conditions are set at the boundary of the upper half-plane:
\beq
V(x_1, 0) = \Phi(x_1), ~ \fr{\prt V}{\prt x_2} (x_1, 0)= \Psi(x_1), ~~\mbox{as} ~ x_1 \in (-\infty, +\infty),
\eequ{eq3a}
where $\Phi$ and $\Psi$ are given bounded functions.

In particular, choose the independent variables $x_1$ and $x_2$ to be a length coordinate $X$ and a time variable $T$, respectively, and refer to the above problem \eq{eq3}--\eq{eq3a} as the Cauchy problem for the wave equation in a temporally and spatially inhomogeneous elastic string. In this case, $V$ stands for the elastic displacement, and the relations \eq{eq3a} represent the initial conditions for the displacement and velocity at $T=0.$ 

\subsection{Temporally inhomogeneous medium}
\label{TIM}

\begin{figure}
\centering
\includegraphics[height=4.0cm]{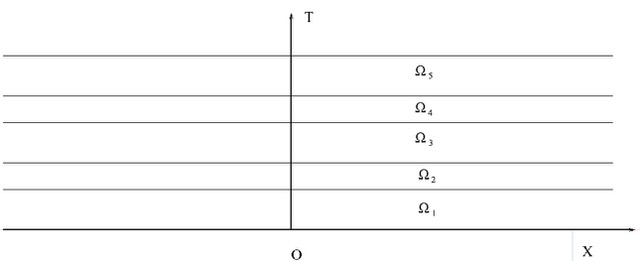}
\caption{\footnotesize Temporal laminate}
\label{Fig1}
\end{figure}

Consider the special case where $\Ga$ and $\Gb$ are $X-$independent, but may change with $T$, i.e. $\Ga=\Ga(T), ~ \Gb = \Gb(T)$. 
In the context of the physical interpretation, related to vibration of an elastic string, which extends along the $X-$axis, the coefficients $\Ga$ and $\Gb$ represent the elastic stiffness and the mass density, respectively. 
Firstly, consider the case when these are continuous and bounded functions of $T$. 
The wave equation  \eq{eq3} then becomes
\beq
\Ga(T) \fr{\prt^2 }{\prt X^2} V(X,T)  - \fr{\prt}{\prt T} \B(\Gb(T) \fr{\prt }{\prt T} V(X,T) \B) = 0.
\eequ{eq4}
It is assumed that the materials are non-dispersive, since otherwise the governing equation \eq{eq4} would take a different form (see, for example, \cite{Engheta}).
Note that in \eq{eq4} by interchanging $T$ and $X$ we obtain the governing equation in a spatial laminate. This allows one to map results for spatial laminates to temporal laminates. The difference is that one has to respect causality: an instantaneous point source can only generate waves in the future, not the past. The response to such point sources, or periodic arrays of them, will be the main focus of our investigations in subsequent sections.

If one has a band gap in the spatial laminate then, in the temporal laminate, one may have, for a fixed wavevector and for a frequency in the bandgap, complex conjugate solutions,
one of which corresponds to waves growing exponentially in time.

The formal application of the Fourier transform 
\beq
\hat{V}(T, k) = \int_{-\infty}^\infty V(X,T) e^{i X k} dX
\eequ{eq4a}
leads to
\beq
\Ga(T) k^2 \hat{V}(k,T)  + \fr{\prt}{\prt T} \B(\Gb(T) \fr{d }{d T} \hat{V}(k,T) \B) = 0.
\eequ{eq5}
In particular, if $\Gb = 1$, and the normalised stiffness coefficient $\Ga$ is a periodic function of time, for example $\Ga(T) = 1 - 2 q \cos (2T), ~ 0 < q < 1/2$, then the equation \eq{eq5} becomes
$$
\fr{d^2 }{d T^2} \hat{V}(k,T)  +  k^2 ( 1 - 2 q \cos (2T)) \hat{V}(k,T)   = 0,
$$  
which is the well-known Mathieu differential equation. The solution can be represented in terms of the Mathieu functions, and Floquet theory can be used accordingly (see, for example, \cite{McLachlan}).  

Consider an alternative configuration, where both the stiffness coefficient $\Ga$ and the mass density coefficient $\Gb$ are piece-wise constant functions of $T$, which implies a ``temporal stratification'' (see Fig. \ref{Fig1}). In this case, the interface conditions should be set at every temporal interface, as discussed in Section  \ref{TS}.

In the general case, equation \eq{eq5} can be reduced to 
\beq
\fr{\prt}{\prt T} \begin{pmatrix} \hat{W} \cr \hat{V} \end{pmatrix} + \begin{pmatrix} 0 & - \tilde{\Ga}(k, T) \cr ({\Gb}(T))^{-1} & 0  \end{pmatrix} \begin{pmatrix} \hat{W} \cr \hat{V} \end{pmatrix} = 0,
\eequ{eq6}
where 
$$
\tilde{\Ga}(\xi,T) = k^2 \Ga(T) ~ \mbox{and} ~ \hat{W}(k, T) = - \Gb(T) \fr{\prt}{\prt T}  \hat{V}(k, T).
$$
The relations \eq{eq3a} yield the initial conditions for $ \hat{W}$ and $\hat{V}$:
\beq
\hat{V}\B|_{T=0} = \hat{\Phi}(k), ~ \hat{W}\B|_{T=0} = \Gb(0)  \hat{\Psi}(k),
\eequ{eq7}
where $\hat{\Phi}(k), \hat{\Psi}(k)$ are the Fourier transforms of the right-hand sides in \eq{eq3a}.

Using the matrix notations
\beq
\BY(T) = \begin{pmatrix} \hat{W} \cr \hat{V} \end{pmatrix}, ~ \BCM(T) = \begin{pmatrix} 0 & - \tilde{\Ga}(k, T) \cr ({\Gb}(T))^{-1} & 0 \end{pmatrix},
\eequ{eq8}
and regarding $k$ as a fixed parameter, we have the following initial boundary value problem
\beq
\fr{d}{d T} \BY (T) = - \BCM(T) \BY(T), 
\eequ{eq9}
\beq
\BY(0) = \BY_0 = \begin{pmatrix} \hat{W} \cr \hat{V} \end{pmatrix} \B|_{T=0}.
\eequ{eq9a}
Hence, the vector function $\BY (T)$ evaluated at $T=T_*$ can be written in the form
\beq
\BY(T_*) = \exp(-\int_0^{T_*} \BCM(\tau) d \tau) \BY_0.
\eequ{eq10}
 In the case where 
$\BCM(\tau)$ is piecewise constant in $\tau$  the exponential term can be written as a product of matrices. This corresponds to the transfer matrix approach in spatial laminates. 

\subsection{Temporal stratification}
\label{TS}

Consider a horizontally stratified half-plane $\BR^2_+= \{(X,T): X \in \BR, ~T >0\} = \overline{\cup_{j=1}^\infty \GO_j}$, as shown in Fig. \ref{Fig1}, where $\GO_j$ are the horizontal non-intersecting layers. Assume that across the horizontal interfaces separating $\GO_n$ and $\GO_{n+1}$ the functions $V$ and $W$ are continuous.

Let us introduce the set of positive constants $T_1, T_2$ in such a way that the straight line $\{(X, T): X \in \BR, T = T_1\}$ represents the interface separating the adjacent strips   $\GO_1$ and $\GO_{2}$, and the straight line $\{(X, T): X \in \BR, T = T_1+T_2\}$ represents the interface separating the adjacent strips   $\GO_2$ and $\GO_{3}$. The stratified structure is assumed to be periodic so that the line $\{(X, T): X \in \BR, T = n(T_1+T_2)\}$ separates the strips $\GO_{2n}$ and $\GO_{2n+1}$, and the line $\{(X, T): X \in \BR, T = nT_1+(n-1)T_2)\}$ separates the strips $\GO_{2n-1}$ and $\GO_{2n}$. Also assume that the functions $\Ga(T)$ and  $\Gb(T)$ are piecewise constant, i.e. $\Ga = \Ga_1$ and  $\Gb = \Gb_1$ when $(X, T) \in \GO_{2 n-1}$, and $\Ga = \Ga_2$ and  $\Gb = \Gb_2$ when $(X, T) \in \GO_{2 n}$ for all $n \in \BN$.

To obtain a field pattern, and simplify the analysis, we consider the case where
\beq
T_1 \sqrt{\Ga_1/\Gb_1} = T_2 \sqrt{\Ga_2/\Gb_2} = d.
\eequ{eq11}
This implies that the wave propagating with the speed of $\sqrt{\Ga_1/\Gb_1}$ covers during the time $T_1$ the same distance $d$ as does the wave, which propagates with the speed $\sqrt{\Ga_2/\Gb_2}$ during the time $T_2$.

In this case, for a strip $\GO_n$  equation \eq{eq10} implies
\beq
\BY(T^{(n)}_+) = \fr{1}{2}\begin{pmatrix}  
e^{id|\xi|} + e^{-id|\xi|} & - i |\xi| \sqrt{\Ga \Gb} (e^{id|\xi|} - e^{-id|\xi|}) \cr
\fr{i}{|\xi| \sqrt{\Ga \Gb}} (e^{id|\xi|} - e^{-id|\xi|}) & e^{id|\xi|} + e^{-id|\xi|}  \end{pmatrix} \BY(T^{(n)}_-) ,
\eequ{eq12}
where $\Ga \Gb$ should be replaced by $\Ga_1 \Gb_1$ for odd $\GO_n$ and by $\Ga_2 \Gb_2$ for even $\GO_n$, respectively; the values
$T^{(n)}_-$ and  $T^{(n)}_+$ stand for the lower and upper limits of $T$ corresponding to the boundaries of the strip $\GO_n$. The matrix in the right-hand side of \eq{eq12} has unit determinant. 

This gives an iterative scheme enabling the Fourier transform of the solution to be found in the temporally stratified half-plane, with wave reflections at the temporal interfaces.

\begin{figure}
\centering
\includegraphics[height=9.0cm]{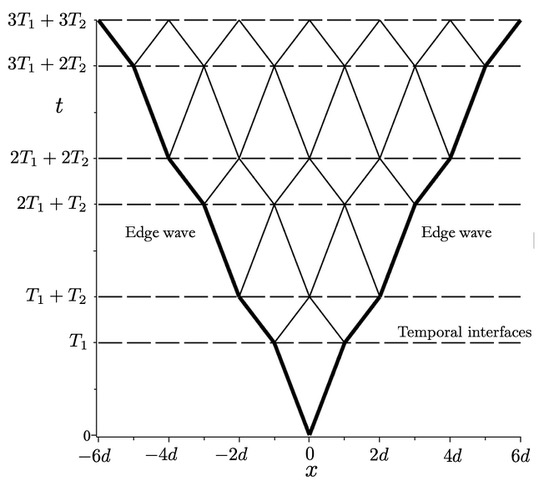}
\caption{\footnotesize Wave splitting at temporal interfaces (dashed lines) showing the emergent field pattern and the two edge waves  (bold lines). The parameter values are $\Ga_1/\Gb_1 =1, \Ga_2/\Gb_2 = 4, T_1=1, T_2 = 1/2.$}
\label{Fig2}
\end{figure}

\section{Field patterns and the edge wave }
\label{edge_wave}

Consider an example demonstrating a  special feature  of the wave, corresponding to the periodic temporal pattern described above.
For the sake of simplicity, assume that $\Psi \equiv 0$ in \eq{eq3a}. In this case, when $(X, T) \in \GO_1$ (i.e. $0 < T < T_1$) the field $V$ is given by
\beq
V(X, T) = \fr{1}{2} \B(\Phi(X- \sqrt{\fr{\Ga_1}{\Gb_1}} T) + \Phi(X+ \sqrt{\fr{\Ga_1}{\Gb_1}} T) \B), ~~ 0 < T < T_1,
\eequ{eq13}

\subsection{The edge wave}
\label{edgewave}

The transmission conditions at $T=T_1$ are continuity of the field $V$ and "linear momentum"  ($\beta \frac{\partial V}{\partial t}$). Application of these conditions leads to the following form of the field $V$ in $\GO_2$, where $T_1 < T < T_1+T_2$ 

\begin{eqnarray}
&V(X, T) = \fr{1}{4} \B\{ (1+\sqrt{\fr{\Ga_1 \Gb_1}{\Ga_2 \Gb_2}}) \B(\Phi (X - \sqrt{\fr{\Ga_2}{\Gb_2}} \tau - d) + \Phi (X + \sqrt{\fr{\Ga_2}{\Gb_2}} \tau + 
d) \B)
&\nonumber \\
&+(1-\sqrt{\fr{\Ga_1 \Gb_1}{\Ga_2 \Gb_2}}) \B(\Phi (X - \sqrt{\fr{\Ga_2}{\Gb_2}} \tau + d) + \Phi (X + \sqrt{\fr{\Ga_2}{\Gb_2}} \tau - d) \B) \B\},&
\label{eq17}
\end{eqnarray}
where the quantity $d$ is defined in \eq{eq11} and $\tau = T - T_1$.

The straightforward observation is that at the temporal interface the wave $V(X,T)$ splits into the ``edge wave'' 
$$
\fr{1}{4} (1+\sqrt{\fr{\Ga_1 \Gb_1}{\Ga_2 \Gb_2}}) \B(\Phi (X - \sqrt{\fr{\Ga_2}{\Gb_2}} \tau - d) + \Phi (X + \sqrt{\fr{\Ga_2}{\Gb_2}} \tau + 
d) \B)
$$
and the ``backward wave''  
$$
\fr{1}{4} (1-\sqrt{\fr{\Ga_1 \Gb_1}{\Ga_2 \Gb_2}}) \B(\Phi (X - \sqrt{\fr{\Ga_2}{\Gb_2}} \tau + d) + \Phi (X + \sqrt{\fr{\Ga_2}{\Gb_2}} \tau - d) \B),
$$
with the coefficient of the ``edge wave'' being greater than the coefficient of the ``backward wave'' within the strip $\GO_2$. We also note that the ``backward wave'' is absent when the coefficients $\Ga, \Gb$ satisfy the relation $\Ga_1 \Gb_1 = \Ga_2 \Gb_2$ corresponding to ``matched impedances". 

\begin{figure}
\centering
\includegraphics[height=6cm]{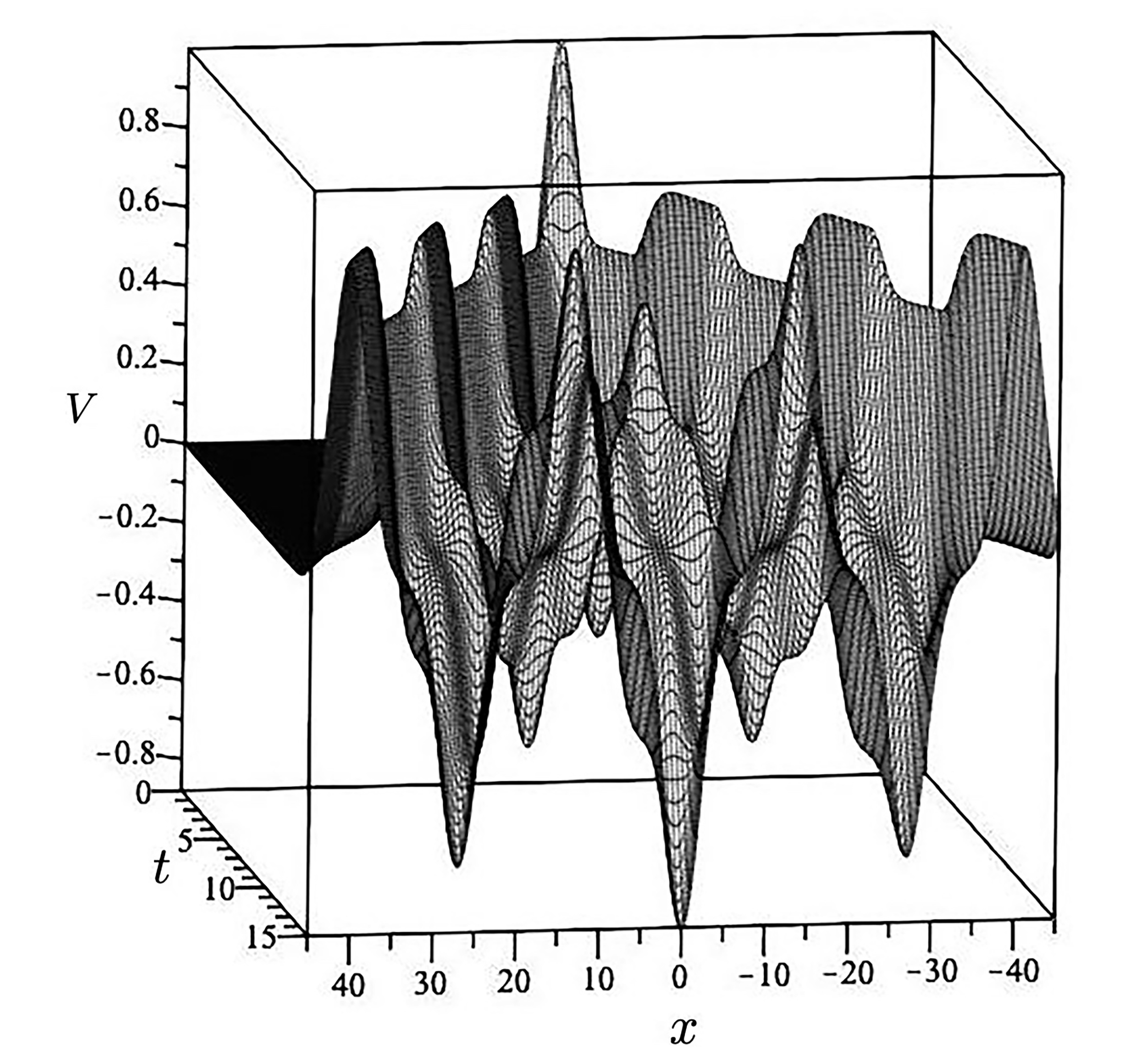} ~~~~~~~~
\includegraphics[height=5.7cm]{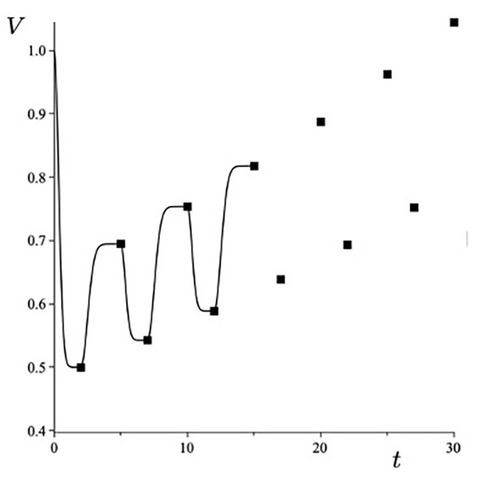}

(a) ~~~~~~~~~~~~~~~~~~~~~~~~~~~~~~~~~~~~~~~~~~~(b) 
\caption{\footnotesize The solution of the Cauchy problem for the case of temporal stratification: 
The following values of the parameters were used in the computations: $\Ga_1 = 8, \Ga_2 = 3, \Gb_1 = 0.7, \Gb_2=0.5906, \kappa = 3.1605, T_1=2, T_2=3.$ The initial profile is given by $\Phi(x)=e^{-x^2/10}$.   (a) The displacement $V(x, t)$, which takes into account both the spatial and temporal dependence. (b) The oscillatory displacement $V$, corresponding to the {\em edge wave};  note the exponential blow-up as $t$ increases: the dots show the representative values of the displacement  in the odd and even laminates, whereas the continuous curve shows the oscillatory behaviour along the edge wave, also shown in the part (a).}
\label{Case0_soln}
\end{figure}

\subsection{Further splitting at the temporal interface}

The temporal interface splitting is illustrated in Fig. \ref{Fig2}, where ray-branching occurs at the times when $\Ga$ and $\Gb$ switch values between $\Ga_1$ and $\Ga_2$, and    $\Gb_1$ and $\Gb_2$, respectively.

The next temporal interface is at $T=T_1+T_2.$ 
By application of the transmission conditions of continuity of the field and "linear momentum" at the temporal interface $T=T_1+T_2$,
the field $V(X,T)$ in the strip $\GO_3$, when $T_1+T_2 < T < 2 T_1 + T_2$, is given by
\begin{eqnarray}
V(X, T) = \fr{1}{8}  \B\{ (\sqrt{\fr{\Ga_2 \Gb_2}{\Ga_1 \Gb_1}} - \sqrt{\fr{\Ga_1 \Gb_1}{\Ga_2 \Gb_2}}) \B[ \Phi \B(X+ \sqrt{\fr{\Ga_1}{\Gb_1}}(T-T_1-T_2)\B) \nonumber \\
+ 
\Phi \B(X-\sqrt{\fr{\Ga_1}{\Gb_1}}(T-T_1-T_2)\B) \B] \nonumber \\
+ (2 + \sqrt{\fr{\Ga_2 \Gb_2}{\Ga_1 \Gb_1}} +  \sqrt{\fr{\Ga_1 \Gb_1}{\Ga_2 \Gb_2}} ) \B[ \Phi \B(X+ \sqrt{\fr{\Ga_1}{\Gb_1}}T + (\sqrt{\fr{\Ga_2}{\Gb_2}} - \sqrt{\fr{\Ga_1}{\Gb_1}})T_2\B) \nonumber \\ + 
\Phi \B(X-\sqrt{\fr{\Ga_1}{\Gb_1}}T - (\sqrt{\fr{\Ga_2}{\Gb_2}} - \sqrt{\fr{\Ga_1}{\Gb_1}})T_2)\B) \B] \nonumber \\
+ (2 - \sqrt{\fr{\Ga_2 \Gb_2}{\Ga_1 \Gb_1}} -  \sqrt{\fr{\Ga_1 \Gb_1}{\Ga_2 \Gb_2}} ) \B[ \Phi \B(X+ \sqrt{\fr{\Ga_1}{\Gb_1}}(T-T_1-T_2) \B) \nonumber \\ + 
\Phi \B(X-\sqrt{\fr{\Ga_1}{\Gb_1}}(T-T_1-T_2)\B) \B] \nonumber \\
-(\sqrt{\fr{\Ga_2 \Gb_2}{\Ga_1 \Gb_1}} - \sqrt{\fr{\Ga_1 \Gb_1}{\Ga_2 \Gb_2}}) \B[ \Phi \B(X- \sqrt{\fr{\Ga_1}{\Gb_1}}(T-T_1-T_2) + 2H \B) \nonumber \\
+ 
\Phi \B(X+\sqrt{\fr{\Ga_1}{\Gb_1}}(T-T_1-T_2)-2H\B) \B] \B\}.
\label{eq20}
\end{eqnarray}
We note that the ``edge wave'' term (in the above formula)
$$
\fr{1}{8} (2 + \sqrt{\fr{\Ga_2 \Gb_2}{\Ga_1 \Gb_1}} +  \sqrt{\fr{\Ga_1 \Gb_1}{\Ga_2 \Gb_2}} ) \B[ \Phi \B(X+ \sqrt{\fr{\Ga_1}{\Gb_1}}T + (\sqrt{\fr{\Ga_2}{\Gb_2}} - \sqrt{\fr{\Ga_1}{\Gb_1}})T_2\B) 
$$
$$ + 
\Phi \B(X-\sqrt{\fr{\Ga_1}{\Gb_1}}T - (\sqrt{\fr{\Ga_2}{\Gb_2}} - \sqrt{\fr{\Ga_1}{\Gb_1}})T_2)\B) \B]
$$
has the largest coefficient
$$
\fr{1}{8} (2 + \sqrt{\fr{\Ga_2 \Gb_2}{\Ga_1 \Gb_1}} +  \sqrt{\fr{\Ga_1 \Gb_1}{\Ga_2 \Gb_2}} )
$$
compared to other coefficients in the representation resulting from the temporal interface split; we also note that some of the coefficients in \eq{eq20} are negative.

%


\subsection{The edge wave blow-up}

The process may be applied to further layers to give a more detailed solution of the Cauchy problem at any given time. The results are shown in Fig.  \ref{Case0_soln}a, where 
the ``edge wave'' is clearly visible.
Although the magnitude of the ``edge wave'' is not monotonic in time,
it may be seen that after passing $n$ macro-cells $\GO_1^{(j)} \cap \GO_2^{(j)}, ~ j=1, \ldots, n$, we have the following amplitude coefficients
\beq
\CC^{(1)}_n = \fr{1}{2} \B( 1 + \fr{1}{4} ( \sqrt{\fr{\Ga_1 \Gb_1}{\Ga_2 \Gb_2}} + \sqrt{\fr{\Ga_2 \Gb_2}{\Ga_1 \Gb_1}} - 2 ) \B)^n 
~~~~\mbox{in} ~~ \GO_1^{(n+1)},\eequ{eq21}
\beq
\CC^{(2)}_n = \fr{1}{2} \Big(1+ \sqrt{\fr{\Ga_1 \Gb_1}{\Ga_2 \Gb_2}} \Big) \CC^{(1)}_n ~~~~ \mbox{in} ~~ \GO_2^{(n+1)},
\eequ{eq21a}
which both grow exponentially, as $n \to \infty$ for all cases where the positive coefficients $\Ga$ and $\Gb$ are chosen in such a way that $\Ga_1 \Gb_1 \neq \Ga_2 \Gb_2$. The formulae \eq{eq21} and \eq{eq21a}, characterising the exponential growth of the {edge wave} amplitude in time, are illustrated in Fig.  \ref{Case0_soln}b.




\section{Field patterns for periodic initial conditions}
\label{per_initial}

Another example involving  periodic initial conditions, leading to possible growth of the solution due to wave splitting at temporal interfaces, is considered below.

Consider initial conditions \eq{eq3a}, with $\Phi = 0$, while $$\Psi(X) = \sum_{k=-\infty}^\infty (-1)^k \Gd (X - kD),$$ with $D$ being a 
positive parameter, which represents half a period of the field pattern along the $X-$axis. Here $\Gd(X)$ is the Dirac delta function.


\subsection{``Reflection'' at temporal interfaces}
\label{1_layer}
For the time interval $(0, T_1)$, it is assumed that  $\Ga= \Ga_1$, $\Gb = \Gb_1$ are constant. Thus, when $0 < T < T_1$ the field $V$ takes the form
\beq
V(X, T) = \fr{1}{2} \sqrt{\fr{\Gb_1}{\Ga_1} } \sum_{k=-\infty}^\infty (-1)^k
\B\{    H(X- k D + \sqrt{\fr{\Ga_1}{\Gb_1}} T) -  H(X- k D - \sqrt{\fr{\Ga_1}{\Gb_1}} T) \B\} .
\eequ{eq22}

The derivation is similar to the one discussed in Section \ref{edge_wave}.
When $ T_1 < T < T_1+T_2$ we assume that $\Ga= \Ga_2$, $\Gb = \Gb_2$ take new constant values, so that  
$$
T_2 \sqrt{\Ga_2 / \Gb_2} = T_1 \sqrt{\Ga_1 / \Gb_1}.
$$ 

Let $\tau=T-T_1$ and consider the field $V(X,T)$  in $\GO_2$, where $T_1 < T < T_1+T_2$ in the form 
\beq
V(X, T) = W_I (X - \sqrt{\fr{\Ga_2}{\Gb_2}} \tau) + W_{II} (X + \sqrt{\fr{\Ga_2}{\Gb_2}} \tau),
\eequ{eq14}
 where the functions $W_I$ and $W_{II}$ represent the waves propagating in the positive and negative direction of the $X-$axis. The  transmission conditions, representing the continuity of the ``displacement'' $V$ and of the ``linear momentum'', are set at the temporal interface $T= T_1 ~ (\tau=0)$:
$$
W_I(X) + W_{II}(X) = \fr{1}{2} \sqrt{\fr{\Gb_1}{\Ga_1}} \sum_{k=-\infty}^\infty (-1)^k 
\B\{    H(X- k D + \sqrt{\fr{\Ga_1}{\Gb_1}} T_1) 
$$
\beq
-  H(X- k D - \sqrt{\fr{\Ga_1}{\Gb_1}} T_1) \B\} ,
\eequ{eq23}
and
$$
\sqrt{\Ga_2 \Gb_2}(W'_{II}(X) - W'_I(X)) = \fr{\Gb_1}{2} \sum_{k=-\infty}^\infty (-1)^k \B\{
\Gd(X- kD + \sqrt{\fr{\Ga_1}{\Gb_1}}T_1) 
$$
\beq
+ \Gd(X- kD - \sqrt{\fr{\Ga_1}{\Gb_1}}T_1) \B\}.
\eequ{eq24}
When $T_1 < T < T_1+T_2$ and $|X| < D/2$, direct integration leads to the following representation of the solution
 $$
 V(X, T) = W_I (X - \sqrt{\Ga_2/\Gb_2} \tau) + W_{II}  (X - \sqrt{\Ga_2/\Gb_2} \tau)
 $$ 
 $$
 = \fr{1}{4} \sqrt{\fr{\Gb_1}{\Ga_1}} \Bigg\{\sum_{k=-\infty}^\infty (-1)^k
 \Big[ H(X -\sqrt{\Ga_2/\Gb_2} \tau - kD +\sqrt{\Ga_1/\Gb_1} T_1)
 $$
 $$
 - H(X -\sqrt{\Ga_2/\Gb_2} \tau - kD -\sqrt{\Ga_1/\Gb_1} T_1)
 +H(X +\sqrt{\Ga_2/\Gb_2} \tau - kD +\sqrt{\Ga_1/\Gb_1} T_1)
$$
$$
-H(X +\sqrt{\Ga_2/\Gb_2} \tau - kD -\sqrt{\Ga_1/\Gb_1} T_1) \Big]
$$
$$
+\sqrt{\fr{\Ga_1 \Gb_1}{\Ga_2 \Gb_2}} \sum_{k=-\infty}^\infty (-1)^k \int_{X-\sqrt{{\Ga_2}/{\Gb_2}}\tau}^{X+\sqrt{{\Ga_2}/{\Gb_2}}\tau} \Bigg( \Gd(\xi - k D + \sqrt{\Ga_1/\Gb_1} T_1) 
$$
\beq
+\Gd(\xi - k D - \sqrt{\Ga_1/\Gb_1} T_1) \Bigg) d \xi \Bigg\} .
\eequ{eq25}

\subsection{The timing of  the ``temporal switch'' versus the spatial periodicity}
\label{2_layer}

We choose $T_1$ and $T_2$ so that 
\beq
\sqrt{\fr{\Ga_1}{\Gb_1}} T_1 = \fr{D}{2} = \sqrt{\fr{\Ga_2}{\Gb_2}} T_2,
\eequ{eq26}
i.e. the first wave reaches the boundary of the elementary cell when the temporal switch occurs. 

Following the derivation, which is similar to that in Section \ref{1_layer}, we deduce that when $T > T_1+T_2$, the field $V(X, T)$ takes the form
$$
V(X, T) = \fr{1}{2} \sqrt{\fr{\Gb_1}{\Ga_1}} \sqrt{ \fr{\Ga_2 \Gb_2}{\Ga_1 \Gb_1}}
\sum_{k=-\infty}^\infty (-1)^k \Bigg(H(X - k D - \sqrt{\Ga_1/\Gb_1} (T-T_2-T_1))
$$
\beq
-H(X - k D - \sqrt{\Ga_1/\Gb_1} (T-T_2-T_1))\Bigg) .
\eequ{eq27}

The above equation is similar to \eq{eq22} subject to the replacement of the argument $T$ in \eq{eq22} by
 $T-T_2-T_1$ in \eq{eq27} and an additional factor $\sqrt{ \fr{\Ga_2 \Gb_2}{\Ga_1 \Gb_1}}$ in \eq{eq27}.
 In particular, when $X=0$ and $T > T_1+T_2$
 $$
 V(0, T) = - \fr{1}{2}\sqrt{\fr{\Gb_1}{\Ga_1}}  \sqrt{ \fr{\Ga_2 \Gb_2}{\Ga_1 \Gb_1}}
 \sum_{k=-\infty}^\infty (-1)^k \Bigg\{ H( - k D - \sqrt{\Ga_1/\Gb_1} (T-T_2-T_1)) 
 $$
 \beq
 - H( k D - \sqrt{\Ga_1/\Gb_1} (T-T_2-T_1)) \Bigg\} = -\sqrt{\fr{\Gb_1}{\Ga_1}} \sqrt{\fr{\Ga_2 \Gb_2}{\Ga_1 \Gb_1}} .
 \eequ{eq28}
 By considering a periodic semi-infinite temporal stratification with the elementary cell of the width $T_1 + T_2$, after $n$ iterations we obtain
 \beq
 V(0, T) = (-1)^n \sqrt{\fr{\Gb_1}{\Ga_1}} \Bigg(\fr{\Ga_2 \Gb_2}{\Ga_1 \Gb_1}\Bigg)^{n/2},
 \eequ{eq29}
 when $n(T_1+T_2) < T < (n+1) T_1 + n T_2.$

The above formula \eq{eq28} shows that when $\fr{\Ga_2 \Gb_2}{\Ga_1 \Gb_1} > 1$, the modulus of the displacement at the origin $|V(0,T)|$ increases exponentially as $n \to \infty.$

  \section{Chiral interfaces}
\label{chiral_interface}

The scalar problem of vibration of an elastic string considered above 
will here be extended to a vector case. Both transverse and longitudinal displacement are included in the formulation, while the coupling process is governed by a chiral term, present in the equation or interface transmission conditions. 

Firstly, the model of a chiral continuum will be summarised since this continuum will be used as a chiral temporal interface. Secondly, energy considerations will be discussed followed by the analysis of a transient configuration of a spatial chiral interface.

\subsection{The chiral medium}
\label{chiral_medium}

A model for a chiral medium has been introduced in stages. Firstly, gyroscopic resonators connected to elastic beams were introduced in \cite{GyroMS} where a gyroscopic resonator was replaced by suitable displacement boundary conditions replacing the resonator and describing the gyroscopic action. The linearised framework of the gyroscopic motion was used and a small angle of nutation was assumed. The governing equations for a one-dimensional discrete chain of such resonators with hinged bases and connected by massless springs was discussed in \cite{Jones2020} and the resonators were again replaced by appropriate time dependent displacement boundary conditions. The system of coupled governing equations in the longitudinal and transverse displacements was homogenised in \cite{Jones2021}  leading to a model with governing equations, written for a chiral medium, in the vector form

\beq
-\fr{\prt^2}{\prt t^2}{\BU} +  \BD \fr{\prt^2}{\prt x^2}\BU + \Ga \BR \fr{\prt }{\prt t}{\BU} = {\bf 0}, 
\eequ{ceq1}

where $\BU = (u(x,t), v(x,t))^T$ is a real-valued vector function, representing the longitudinal, $u(x,t)$, and  transverse, $v(x,t)$, displacements, 
\beq
\BD = \mbox{diag} \{c_1^2, c_2^2   \}, ~~ \BR = \begin{pmatrix} 0 & 1 \cr -1 & 0\end{pmatrix},
\eequ{ceq2a}
where $c_1$ and $c_2$ are the wave speeds for the longitudinal and transverse waves respectively, and $\Ga$ is the (re-defined) gyricity parameter. The third term in  \eq{ceq1} is the coupling term, and it represents the gyroscopic force, which is orthogonal to the velocity vector.

Although, in the linearised setting, the gyroscopic force appears to be non-conservative, it is orthogonal to the velocity vector, and for a finite chiral elastic rod the classical energy conservation holds, as demonstrated in the illustrative example below.

\subsection{Energy consideration} 

Here, an elementary demonstration is given, based on two examples which include vibration of a chiral finite rod and time-harmonic vibration of a single chiral inertial resonator. 

\subsubsection{The finite elastic chiral rod}

Let $\BU = (u(x,t), v(x,t))^T$ be a real-valued vector function, representing the displacements, which satisfy the equations of motion (\ref{ceq1}) for the homogenised chiral elastic rod, $x \in (0,1)$
with the boundary conditions
\beq
\BU(0, t) = \BU(1,t) = 0 ~~\mbox{for any admissible } t.
\eequ{ceq2}

Multiplying equation \eq{ceq1} by $\fr{\prt}{\prt t} \BU^T$ and integrating with respect to $x$ over the interval $(0, 1)$ gives
\beq
0 = \int_0^1 \Big\{  \fr{\prt u}{\prt t}     \fr{\prt^2 u}{\prt t^2}  + \fr{\prt v}{\prt t}     \fr{\prt^2 v}{\prt t^2}  -  c_1^2  \fr{\prt u}{\prt t}     \fr{\prt^2 u}{\prt x^2}  - c_2^2 \fr{\prt v}{\prt t}     \fr{\prt^2 v}{\prt x^2}     \Big\} d x.
\eequ{ceq3}
Note that $  (\fr{\prt }{\prt t}{\BU})^T \BR \fr{\prt }{\prt t}{\BU} = {\bf 0}.$ By integrating \eq{ceq3} by parts, and using the boundary conditions \eq{ceq2}, we obtain
\beq
\fr{d}{d t} \Big( K + P   \Big) = 0,
\eequ{ceq4}
where $K$ and $P$ represent the kinetic and potential energies, respectively, 
\beq
K = \fr{1}{2} \int_0^1 \Big( (\fr{\prt u}{\prt t})^2 + (\fr{\prt v}{\prt t})^2    \Big) dx,  ~ P = \fr{1}{2} \int_0^1 \Big(  c_1^2 (\fr{\prt u}{\prt x})^2 + c_2^2 (\fr{\prt v}{\prt x})^2    \Big) dx.
\eequ{ceq5}
Hence, the standard conservation law holds, and the total energy $K + P$ is time-independent.

\subsubsection{Time-harmonic regime}

In the time-harmonic regime, with radian frequency $\Go$, the energy consideration for a finite elastic rod leads to the evaluation of the first eigenvalue. Assuming that
$\BU(x, t) = \hat{\BU}(x) \exp(-i \Go t),$ then
\beq
\Go^2 \hat{\BU}(x) + \BD \hat{\BU}''(x) -i\Go \Ga \BR \hat{\BU}(x) = {\bf 0}, ~ x \in (0, 1),
\eequ{ceq6}
\beq
\hat{\BU}(0) = \hat{\BU}(1) = {\bf 0}.
\eequ{ceq7} 
In the above equations, $\hat{\BU}(x) = (\hat{u}(x), \hat{v}(x))^T$ is a complex valued vector function.
When the gyricity parameter $\Ga$ is zero, the problem \eq{ceq6}, \eq{ceq7} splits into two uncoupled standard eigenvalue problems for two harmonic oscillators. 
Additionally, in the case when $\Ga = 0,$ the following identity holds:
\beq
\Go = \Bigg(\fr{\int_0^1 \Big( c_1^2 |\hat{u}'(x)|^2 +   c_2^2 |\hat{v}'(x)|^2   \Big) dx}{\int_0^1 \Big( |\hat{u}(x)|^2 +   |\hat{v}(x)|^2   \Big) dx}\Bigg)^{1/2},
\eequ{ceq8}
where the energy integral is in the numerator of \eq{ceq8}. Using the normalisation  
\beq 
\int_0^1 \Big( |\hat{u}(x)|^2 +   |\hat{v}(x)|^2   \Big) dx =1,
\eequ{ceq8a}
and adopting the notation $H_0(0,1)$ for the space of  vector functions, which satisfy \eq{ceq7}, \eq{ceq8a} and have the finite energy integral in \eq{ceq8}, we have that the first eigenvalue as
\beq
\Go_1 = \inf_
{\hat{\BU} \in H_0(0,1)}  \Bigg({\int_0^1 \Big( c_1^2 |\hat{u}'(x)|^2 +   c_2^2 |\hat{v}'(x)|^2   \Big) dx} \Bigg)^{1/2}.
\eequ{ceq9}
When $\Ga \neq 0 $ and the normalisation \eq{ceq8a} is in place, for positive $\Go$ the identity \eq{ceq8} is replaced by
\beq
\Go = \fr{\int_0^1 \Big( c_1^2 |\hat{u}'(x)|^2 +   c_2^2 |\hat{v}'(x)|^2   \Big) dx}{ \Ga \int_0^1 \Im(\bar{\hat{u}}\hat{v}) dx +
\sqrt{\Ga^2  \Big(\int_0^1 \Im(\bar{\hat{u}}\hat{v}) dx \Big)^2   + \int_0^1 \Big( c_1^2 |\hat{u}'(x)|^2 +   c_2^2 |\hat{v}'(x)|^2   \Big) dx } } .
\eequ{ceq10}

\subsection{Transient scattering on an active chiral spatial interface}

 Consider two one-dimensional semi-infinite elastic rods joined by a one-dimensional chiral 
segment. The uncoupled longitudinal and transverse displacements in the rods are governed by (\ref{ceq1}) with $\alpha=0$. The coupled displacements in the chiral 
segment are governed by (\ref{ceq1}) with the chirality parameter $\alpha$ being non-zero.  



For simplicity, the limiting case when the chiral 
segment height tends to zero will be considered, with $\alpha$ approaching infinity in such a way that the limit is an imperfect interface.

\subsubsection{Imperfect spatial chiral interface}
\label{impchsp}

Assume that the chiral interface extends over the small interval $0 < x < d$, $ d \ll 1$, and the following equations and interface transmission conditions hold 
\beq
\fr{\prt^2}{\prt t^2} \BU(x, t) - \Ga \BR \fr{\prt}{\prt t} \BU(x, t) - \BD \fr{\prt^2}{\prt x^2} \BU = {\bf 0}, ~~ \mbox{when}~ 0<x< d,  
\eequ{ceq11a}
\beq
\fr{\prt^2}{\prt t^2} \BU(x, t) - \BD \fr{\prt^2}{\prt x^2} \BU = {\bf 0}, ~~ \mbox{outside the segment} ~ [0, d],
\eequ{ceq11b}
and the continuity transmission conditions are set at $x=0$ and $x=d$
\beq
\Big[ \BU   \Big]_{x=-0}^{x=+0} = {\bf 0}, ~ \Big[ \fr{\prt}{\prt x} \BU   \Big]_{x=-0}^{x=+0} = {\bf 0},
\eequ{ceq11c}
\beq
\Big[ \BU   \Big]_{x=d-0}^{x=d+0} = {\bf 0}, ~ \Big[ \fr{\prt}{\prt x} \BU   \Big]_{x=d-0}^{x=d+0} = {\bf 0}.
\eequ{ceq11d}
Assuming that the gyricity parameter $\Ga$ is large, so that the product $\Gv := \Ga d$ remains constant, and introducing the scaled variable $\xi = x/d$, we can re-write the equation \eq{ceq11a} in the form:
\beq
\BD \fr{\prt^2}{\prt \xi^2} \BU + d \Gb \BR \fr{\prt}{\prt t} \BU - d^2 \fr{\prt^2}{\prt t^2} \BU = {\bf 0}, ~~ 0< \xi < 1.  
\eequ{ceq11e}
Consider the asymptotic approximations of the form
\beq
\BU = \BU^{(0)} (\xi, t) + d \BU^{(1)} (\xi, t) + O(d^2), ~~\mbox{when}~ 0 < \xi < 1,
\eequ{ceq11f}
and
\beq
\BU = \BV^{(0)}(x, t) + O(d), ~~ \mbox{when} ~ x < 0 ~\mbox{or} ~x > d. 
\eequ{ceq11g}
Substitution of \eq{ceq11f} into \eq{ceq11a}, and use of the transmission conditions \eq{ceq11c} and \eq{ceq11d} leads to a sequence of problems for $\BU^{(0)}$ and $\BU^{(1)}$. Namely, we deduce
 \beq
\fr{\prt^2}{\prt \xi^2} \BU^{(0)} = {\bf 0}, ~~ 0 < \xi < 1,
\eequ{ceq11h}
\beq
\fr{\prt}{\prt \xi} \BU^{(0)} = {\bf 0} ~~\mbox{when} ~ \xi =+0 ~ \mbox{and} ~ \xi = 1-0.
\eequ{ceq11i}
Hence, $\BU^{(0)}$ is $\xi-$independent, and 
\beq
\BU^{(0)}(t) = V^{(0)} (0, t) = V^{(0)}(d, t).
\eequ{ceq11j}
Next, the vector function $\BU^{(1)}$ satisfies the equation
\beq
\BD \fr{\prt^2}{\prt \xi^2} \BU^{(1)}(\xi, t) + \Gb \BR \fr{\prt}{\prt t} \BU^{(0)}(t) = {\bf 0},
\eequ{ceq11k}
and
\beq
\fr{\prt}{\prt \xi} \BU^{(1)}(0, t) = \fr{\prt}{\prt x} \BV^{(0)}(0,t), ~~\fr{\prt}{\prt \xi} \BU^{(1)}(1, t) = \fr{\prt}{\prt x} \BV^{(0)}(d,t).
\eequ{ceq11kl}
When $d$ is infinitesimally small, \eq{ceq11i}-\eq{ceq11kl} lead to the following transmission conditions for the leading-order term $\BV^{(0)}$ in \eq{ceq11g}, across the chiral interface
\beq
\Big[ \BV^{(0)}(x, t) \Big]_{x = -0}^{x=+0} = {\bf 0}, ~ \BD \Big[ \fr{\prt}{\prt x} \BV^{(0)}(x, t)\Big]_{x = -0}^{x=+0}  = -\Gb \BR \fr{\prt}{\prt t} \BV^{(0)}(0,t).
\eequ{ceq11m}
Taking into account the above asymptotic approximation, in the limit, as $d \to +0$, one can consider the transmission problem for an imperfect spatial  chiral interface, where tractions become discontinuous, and a coupling is observed  between the longitudinal and transverse displacements:
\begin{eqnarray}
\BD \fr{\prt^2}{\prt x^2} \BU -\fr{\prt^2}{\prt t^2} \BU(x, t)= {\bf 0}, ~~ \mbox{when} ~ x \neq 0,
\l{ceq11b1} \\
\Big[ \BU(x, t) \Big]_{x = -0}^{x=+0} = {\bf 0}, ~ \BD \Big[ \fr{\prt}{\prt x} \BU(x, t)\Big]_{x = -0}^{x=+0}  = -\Gb \BR \fr{\prt}{\prt t} \BU(0,t).
\l{ceq11b2}
\end{eqnarray}

We also note that equivalently, instead of \eq{ceq11b1} and \eq{ceq11b2}, we can write the equation with a delta-function term as follows:
\beq
\BD \fr{\prt^2}{\prt x^2} \BU  - \fr{\prt^2}{\prt t^2} \BU(x, t) + \Gb \BR \fr{\prt}{\prt t} \BU(0, t) \delta(x) = {\bf 0}.
\eequ{ceq11n}

\subsubsection{A model scattering problem}
Assume that a step-like front of the longitudinal wave is propagating in the positive direction of the $x-$axis
\beq
\BU_{\mbox{\tiny inc}} = \begin{pmatrix}  1-H(x-x_0-c_1 (t-t_0)) \cr   0 \end{pmatrix} ,
\eequ{ceq11}
and at time $t=t_0 >0$ this front meets a chiral interface, positioned at the point $x=x_0,$ and characterised by gyricity $\Gb$, as discussed in Section \ref{impchsp}.

The total field includes two terms
\beq
\BU(x, t) = \BU_{\mbox{\tiny inc}} + \BU_{\mbox{\tiny sc}},
\eequ{ceq12}
where the incident field has only the first non-zero component $u_{\mbox{\tiny inc}}$, which satisfies the homogeneous wave equation
\beq
\fr{\prt^2}{\prt t^2}{u_{\mbox{\tiny inc}}} -  c_1^2 \fr{\prt^2}{\prt x^2} u_{\mbox{\tiny inc}} = {0}, 
\eequ{ceq13a}
whereas the scattered field has only the second non-zero component $v_{\mbox{\tiny sc}}$, corresponding to the transverse vibration, and it is triggered by an instantaneous point force, whose magnitude is proportional to the velocity of the longitudinal vibration but the orientation of the force  is orthogonal to the longitudinal velocity. Taking into account \eq{ceq11n}, we observe that the scattered field is the transverse transient wave, which satisfies the following initial value problem
\beq
  c_2^2 \fr{\prt^2}{\prt x^2} v_{\mbox{\tiny sc}}-\fr{\prt^2}{\prt t^2}{v_{\mbox{\tiny sc}}}  - \Gb \fr{\prt u_{\mbox{\tiny inc}}}{\prt t}(x_0, t) \Gd(x-x_0) = {0}, 
\eequ{ceq13b}
subject to the initial conditions
\beq
v_{\mbox{\tiny sc}} (x, 0) = 0, ~ \fr{\prt}{\prt t} v_{\mbox{\tiny sc}} (x, 0) = 0,
\eequ{ceq13c}
with $\Ga$ being the chirality parameter characterising the chiral point scatterer.
Taking into account \eq{ceq11n} and \eq{ceq11}, the equation \eq{ceq13b} can be re-written in the form
\beq
 c_2^2 \fr{\prt^2}{\prt x^2} v_{\mbox{\tiny sc}}-\fr{\prt^2}{\prt t^2}{v_{\mbox{\tiny sc}}}  - \Gb \Gd(t-t_0) \Gd(x-x_0) = {0}. 
\eequ{ceq13d}
Taking the Fourier transform $\hat{v}_{\mbox{\tiny sc}} (k, t) = \int_{-\infty}^\infty  v_{\mbox{\tiny sc}} (x, t) \exp(i k x) dx $, we deduce
\beq
 k^2 c_2^2  \hat{v}_{\mbox{\tiny sc}} (k, t)+\fr{\prt^2}{\prt t^2}{\hat{v}_{\mbox{\tiny sc}}} (k, t)  + \Gb \Gd(t-t_0) e^{i k x_0} = {0}, 
\eequ{ceq13e}
\beq
\hat{v}_{\mbox{\tiny sc}} (k, 0) = 0, ~ \fr{\prt}{\prt t} \hat{v}_{\mbox{\tiny sc}} (k, 0) = 0,
\eequ{ceq13cc}
and hence
\beq
\hat{v}_{\mbox{\tiny sc}} (k, t) = - \Gb \fr{\sin(k c_2 (t-t_0))}{k c_2} 
H(t - t_0) e^{i k x_0}.
\eequ{ceq13f}
Using the identity
\beq
\int_0^\infty \fr{\sin(k q)}{k} d k = \fr{\pi}{2} \sign(q),
\eequ{ceq13g}
together with the inverse Fourier transform, we obtain
\beq
v_{\mbox{\tiny sc}} (x,t) = \fr{1}{2 \pi} \int_{-\infty}^\infty \hat{v}_{\mbox{\tiny sc}} (k, t)  e^{-i k x} d k 
\eequ{ceq13h}
$$
= -\fr{\Gb}{2 \pi c_2} H(t-t_0) \int_{-\infty}^\infty \fr{\sin(k c_2 (t-t_0))}{k} e^{i k (x_0-x)} d k
$$
$$
= -\fr{\Gb}{4 c_2} H(t-t_0) \Big(\sign(c_2 (t-t_0) + x - x_0) +  \sign(c_2 (t-t_0) + x_0 - x) \Big) 
$$
$$
= -\fr{\Gb}{2 c_2} H(c_2 (t-t_0) - |x-x_0|)
$$ 
$$
= -\fr{\Gb}{2 c_2} \Big(  H(x-x_0 + c_2 (t-t_0)) -   H(x-x_0 - c_2 (t-t_0))  \Big)
$$

The resulting scattered field represents the transverse wave, as follows
\beq
\BU_{\mbox{\tiny sc}}(x, t) = \begin{pmatrix}  0 \cr  \fr{\Gb }{2 c_2} \Big( H(x-x_0-c_2 (t-t_0)) - H(x-x_0+c_2 (t-t_0)) \Big) \end{pmatrix}.
\eequ{ceq14}
 
We note that the above scattered field is independent of $c_1$. This is not a general result, but a special feature of the model problem considered here, where the incident field is defined via a step function and the homogeneity properties $H(\Gl x) = H(x), \Gd(\Gl x) =  \Gl^{-1} \Gd(x), \Gl > 0,$ have been used.

In the general case of an incident wave, the scattered field may depend on both $c_1$ and $c_2$.

Such a chiral interface can be viewed as an active interface, as it initiates a force, orthogonal to the orientation of the velocity of the incident field, and the magnitude of this force depends on the gyricity parameter $\Gb$. In this case, the term ``active interface'' is used to emphasise that additional energy may enter the system as a result of the interaction of the incident wave with the chiral scatterer. 

For $x < x_0,$ the resulting reflected field $\BU_{\mbox{\tiny refl}}$ has a different polarisation compared to the incident field, and it propagates with  speed $c_2$ compared to the speed $c_1$ of the incident wave
\beq
\BU_{\mbox{\tiny refl}}  (x, t) = \BU_{\mbox{\tiny sc}}  (x, t) = \begin{pmatrix}  0 \cr  -\fr{\Gb}{2 c_2}  H(x-x_0+c_2(t-t_0))  \end{pmatrix}, ~ x < x_0.
\eequ{ceq15}
On the other hand, for $x > x_0$, the transmitted wave has both components, representing longitudinal and transverse waves, as follows
$$
\BU_{\mbox{\tiny transm}}  (x, t) = \BU_{\mbox{\tiny inc}}  (x, t) + \BU_{\mbox{\tiny sc}}  (x, t)   
$$
\beq
 =  \begin{pmatrix}  1-H(x-x_0-c_1 (t-t_0)) \cr  \fr{\Gb }{2 c_2}  (H(x-x_0-c_2(t-t_0))-1)  \end{pmatrix}, ~ x > x_0.
\eequ{ceq16}

We note that the coupling between the longitudinal and transverse vibrations is an essential feature of the chiral scatterer. In the above example, for large values of $\Gb$ 
the chiral transient scatterer can be interpreted as a ``switch'', which initiates a transverse wave, with  amplitude controlled by gyricity,  that may be significantly larger than that of the incident longitudinal wave.  

\section{Scattering on a temporal chiral interface}
\label{temporal_chiral}

In this section, the reflection of a wave on a temporal chiral interface is considered.  

One can re-write the equations of motion \eq{ceq1} in the form.
\beq
\fr{\prt }{\prt t} \Big( \BCM (\Ga t) \fr{\prt}{\prt t} \BU(x, t) \Big) - \BCM(\Ga t) \BD \fr{\prt^2}{\prt x^2} \BU(x,t) = {\bf 0},
\eequ{ceq17}
where the matrix function $\BCM(\Ga t)$ is the rotation matrix defined by 
\beq
\BCM(\Ga t) = \begin{pmatrix}  \cos(\Ga t) & - \sin(\Ga t) \cr \sin(\Ga t) & \cos(\Ga t) \end{pmatrix} .
\eequ{ceq18}

By considering a transient process, we assume that the gyricity parameter $\Ga$ may change instantaneously across a temporal interface $t = T_*$, while the displacement and the momentum remain continuous:
\beq
\Big[  \BU \Big]_{t = T_*-0}^{t=T_*+0} = {\bf 0}, ~~  \Big[ \BCM \fr{\prt}{\prt t}  \BU \Big]_{t = T_*-0}^{t=T_*+0} = {\bf 0}.
\eequ{ceq19} 

\subsection{Transient solutions for small and large values of the gyricity parameter $\Ga$}

Here, illustrations will be presented for solutions of the Cauchy problems on the infinite axis $x \in \BR$ and $t > 0,$ corresponding to different chiral regimes.

The role of the chirality is in the rotational coupling between the longitudinal and transverse displacements and this leads to  wave dispersion. 
In particular if $\Ga=0$, then the problem is reduced to a system of uncoupled Cauchy problems, governed by standard D'Alembert's solutions, which describe non-dispersive waves propagating with the speeds $c_1$ and $c_2$ for the longitudinal and transverse vibrations respectively.

On the other hand, when $\Ga \gg 1$, the solution becomes highly oscillatory in time $t$, compared to the spatial variation in $x$. It can be observed, that for a finite time interval the $x-$dependence appears to be ``frozen''. On a given cross-section, the variable $x$ can be considered as a fixed parameter. 


For the purposes of numerical illustration, it is convenient to express the governing equations (\ref{ceq1}) in dimensionless form. Firstly, the following initial conditions will be assumed:

\begin{equation}
u(x,0)=Lf_1(\frac{x}{L}),~~v(x,0)=Lf_2(\frac{x}{L}),~~
\label{ics}
\end{equation}
$$
\frac{\partial u}{\partial t}(x,0)=c_1 g_1(\frac{x}{L}),~~\frac{\partial v}{\partial t}(x,0)=c_1 g_2(\frac{x}{L}).
$$

Here 
$L$ is defined as a characteristic length introduced in the initial conditions. 
The functions $f_i$ and $g_i$ are dimensionless. 

Introduce the dimensionless variables:   
\begin{equation}
\tilde{x}= \frac{x}{L},~~\tilde{t}= \frac{c_1}{L} t,~~\tilde{u}= \frac{u}{L},~~ \tilde{v}= \frac{v}{L} ,~~\lambda= \frac{c_2}{c_1},~~\gamma = \frac{\alpha L}{c_1} . 
\label{scale}
\end{equation}

The governing equations then become (dropping the tildes $\tilde{~}$)

\begin{equation}
\begin{aligned}
\frac{\partial ^2 u}{\partial t^2}- \frac{\partial ^2 u}{\partial x^2} -\gamma \frac{\partial v}{\partial t}&= 0,  \\
\frac{\partial ^2 v}{\partial t^2}- \lambda^2 \frac{\partial ^2 v}{\partial x^2} +\gamma \frac{\partial u}{\partial t}&= 0,
\end{aligned}\label{dimgoveqn}
\end{equation} 

with initial conditions

\begin{equation}
u(x,0)=f_1(x) 
,~~v(x,0)=f_2(x) 
,~~\frac{\partial u}{\partial t}(x,0)= g_1(x) 
,~~\frac{\partial v}{\partial t}(x,0)= g_2(x). 
\label{dimics}
\end{equation}

It is noted that the system \eq{dimgoveqn} can be written in the matrix form, similar to \eq{ceq17}:
\beq
\fr{\prt}{\prt t} \B( \BCM (\Gg t) \fr{\prt}{\prt t} \BU(x, t) \B)  - \BCM(\Gg t) \BD_\Gl \fr{\prt^2}{\prt x^2} \BU(x, t) = {\bf 0},
\eequ{ceq17a}
where the matrix function $\BCM$ is defined by \eq{ceq18}, and $\BD_\Gl = \mbox{diag} ( 1, \Gl^2 ).$ In this case, it is assumed that $0 < \Gl < 1.$

As discussed in \cite{Jones2020}, for large values of the gyricity, this problem is a singularly perturbed one and an alternative normalisation is used in \cite{Jones2021}, which highlights the highly oscillatory behaviour  of the solutions in time. 

\begin{figure}[H]
\centering
\begin{subfigure}{.5\textwidth}
  \centering
  \includegraphics[width=.9\linewidth]{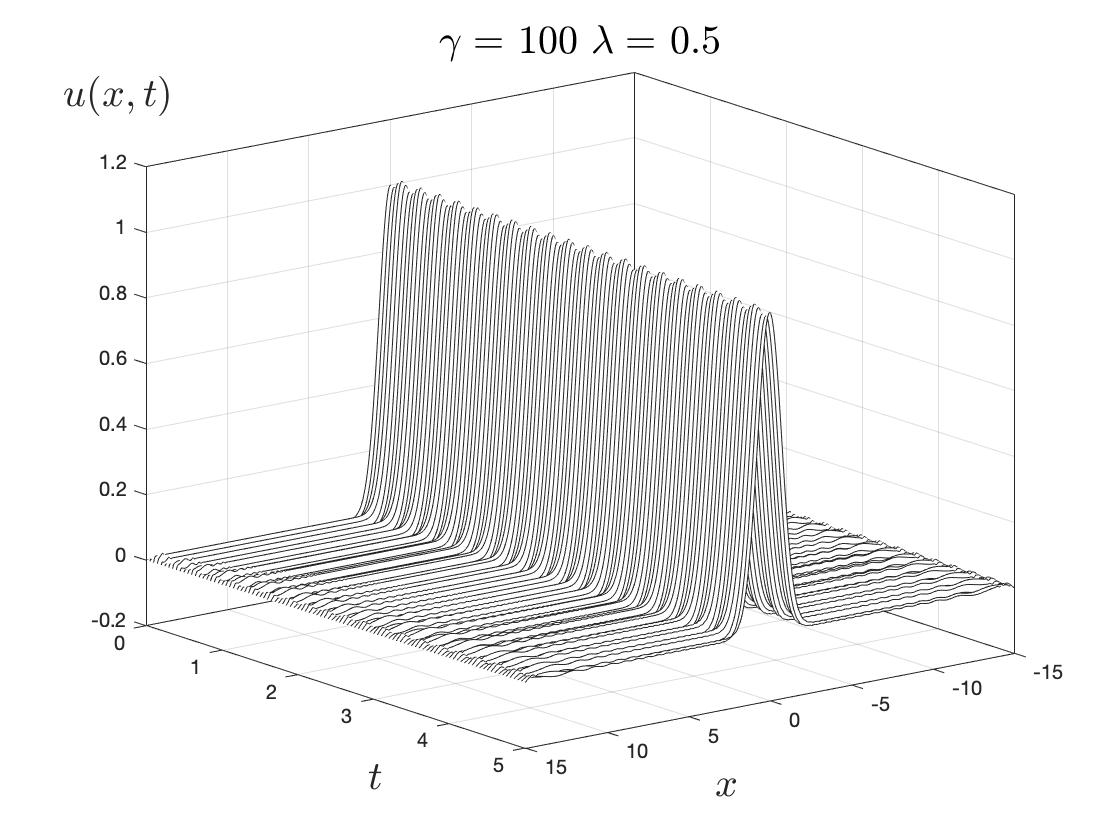}
 
  \caption{\footnotesize The longitudinal displacement $u(x,t)$ .}
  \label{fig1a}
\end{subfigure}%
\begin{subfigure}{.5\textwidth}
  \centering
  \includegraphics[width=.9\linewidth]{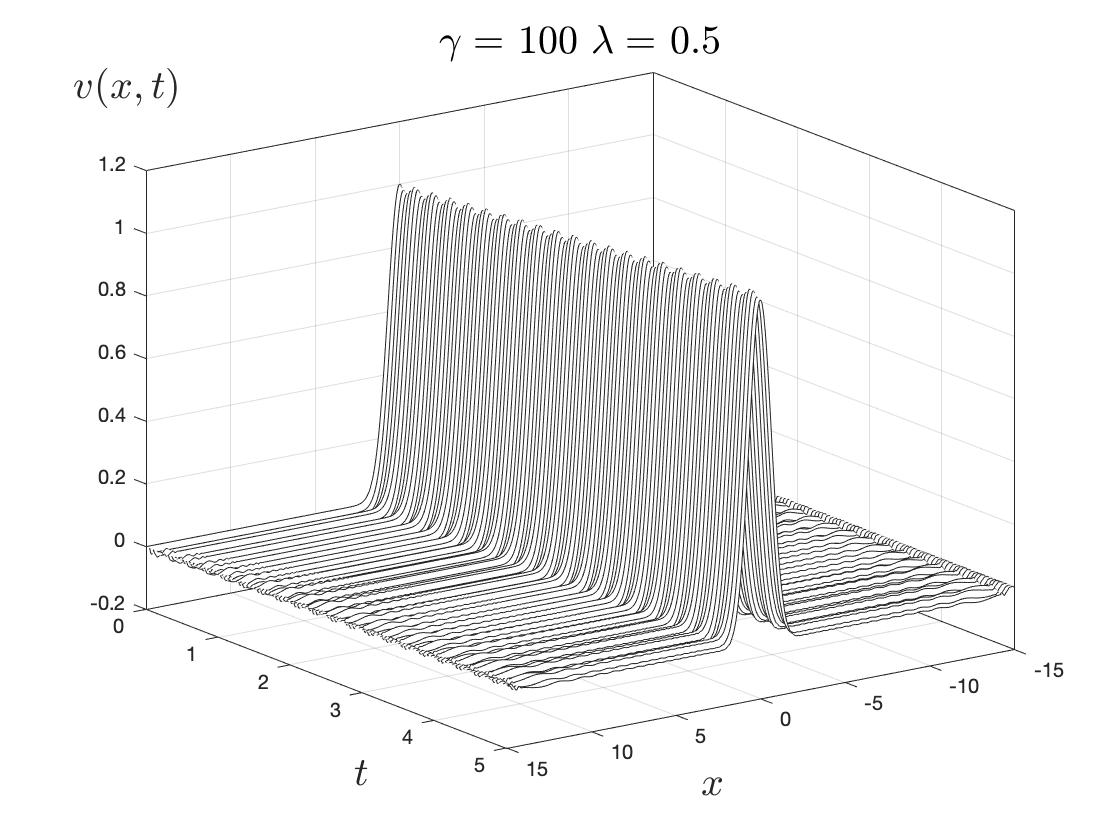}
  \caption{\footnotesize The transverse displacement $v(x,t)$ .}
  \label{fig1b}
\end{subfigure}
\newline
\newline

\caption{\footnotesize Displacements for the case of very large chirality ($\gamma=100$) near the 'stationary limit' between $t=0$ and $t=10$. The wave speed parameter is given by $\lambda=0.5$.}

\label{fig1}
\end{figure}

The "stationary limit" of very large chirality is illustrated by the results shown in Fig. \ref{fig1}. Equations (\ref{dimgoveqn}) have been solved numerically and the dimensionless initial conditions are set as

\begin{equation}
u(x,0)=v(x,0)=e^{-x^2},~~\frac{\partial u}{\partial t}(x,0)=1,~~\frac{\partial v}{\partial t}(x,0)=1 ,
\label{ics2}
\end{equation}

with $\gamma=100$ and $\lambda=0.5$. A dimensionless time frame up to $t=5$ is chosen and an $x$ domain is chosen to be sufficiently large  to minimise reflections. The dimensionless longitudinal and transverse displacements are shown in Figs. \ref{fig1a} and  \ref{fig1b}, respectively. At this relatively large value of $\gamma$, it may be seen that the waves remain approximately stationary in the spatial coordinate whilst they vary harmonically in time. The analytical solution for an infinitely large value of $\alpha$ is given in \cite{Jones2021} and discussed later in this paper.

\begin{figure}[H]
\centering
\begin{subfigure}{.5\textwidth}
  \centering
  \includegraphics[width=.9\linewidth]{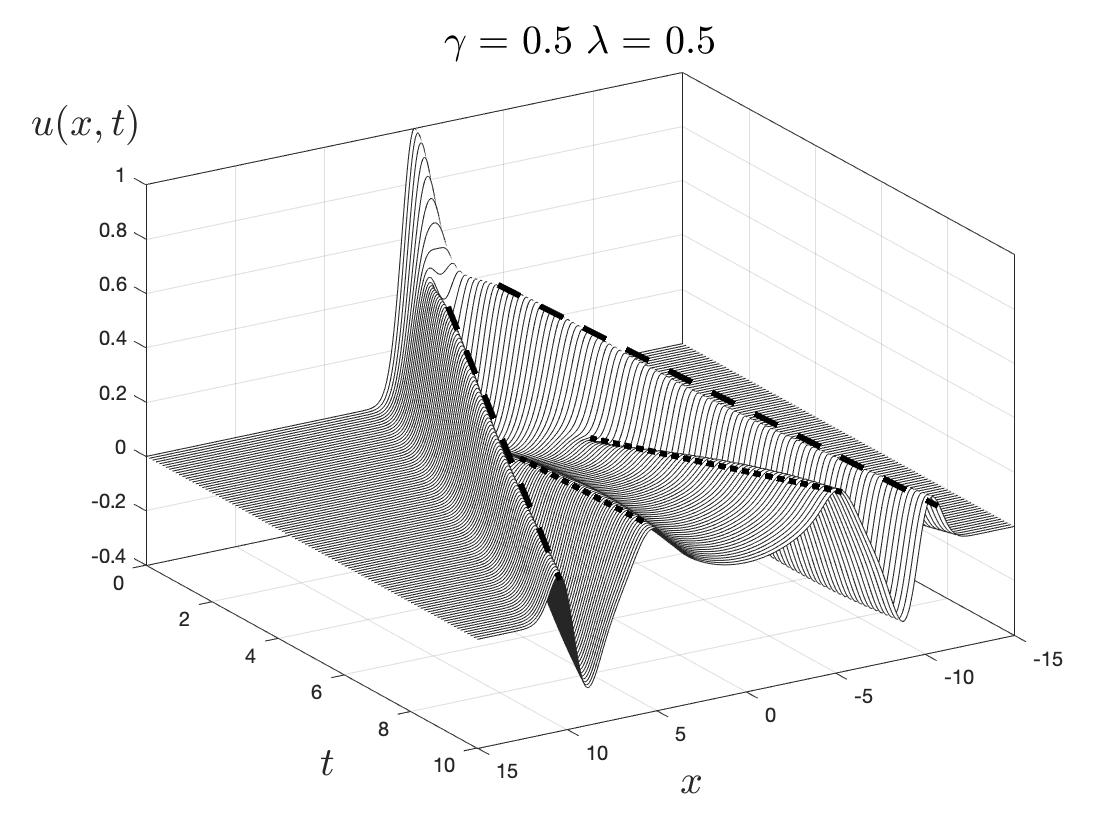}
  \caption{\footnotesize The longitudinal displacement $u(x,t)$.}
  \label{fig2a}
\end{subfigure}%
\begin{subfigure}{.5\textwidth}
  \centering
  \includegraphics[width=.9\linewidth]{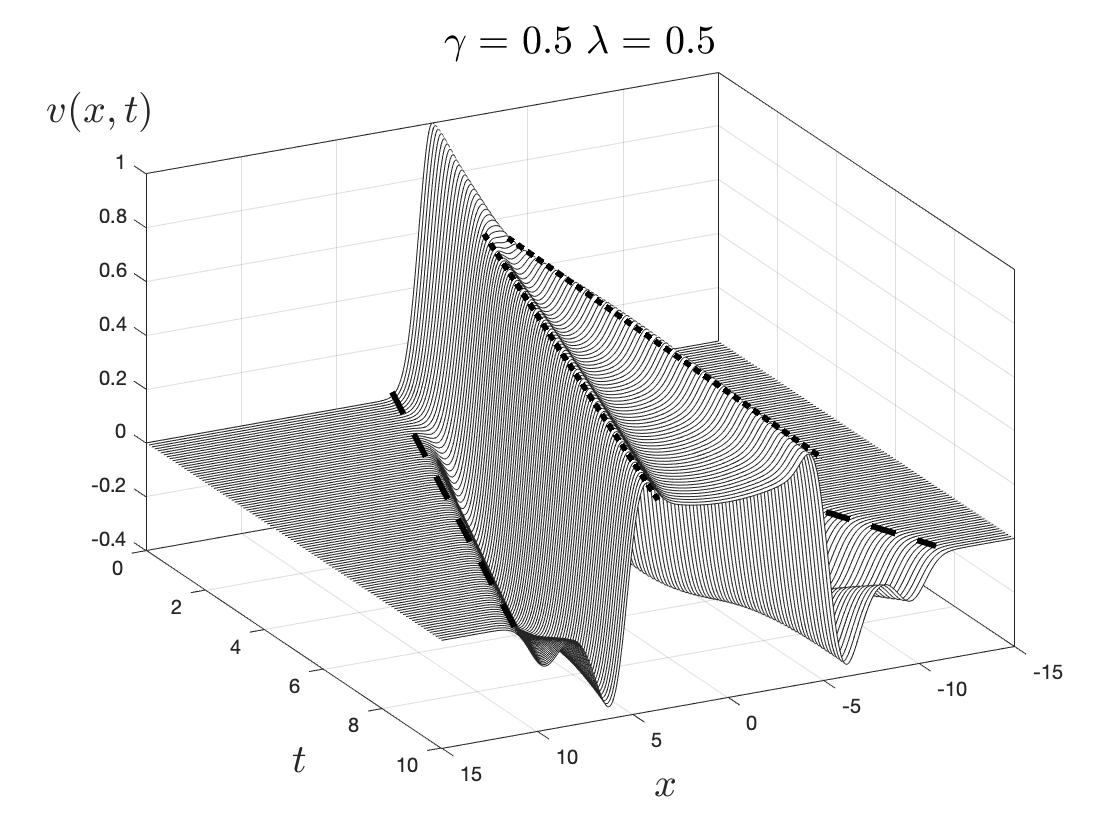}
  \centering
  \caption{\footnotesize The transverse displacement $v(x,t)$.}
  \label{fig2b}
\end{subfigure}
\caption{\footnotesize Longitudinal and transverse displacements between $t=0$ and $t=10$. The relative wave speed parameter $\lambda=0.5$, and the chirality parameter $\gamma=0.5$. The non-chiral characteristics are shown in dashed and dotted lines.}
\label{fig2}
\end{figure}

An example of the case of very small chirality is shown in Fig.  \ref{fig2}. The initial conditions for the displacements are given in (\ref{ics2}) but the initial velocities are set to zero. The relative chirality parameter is given as $\gamma=0.5$ and again  
$\lambda=0.5$. The longitudinal and transverse displacements are shown in Figs. \ref{fig2a} and \ref{fig2b} respectively as functions of $x$ and $t$. It is apparent that the initial wave profile moves to the left and right with increasing time but the shape of that initial profile changes due to the dispersive nature of the waves.

This situation of small chirality is much closer to that of two decoupled wave equations. For the latter case, two Gaussian profiles move to the left and right for both the longitudinal and transverse displacements without dispersion. The respective wave speeds in this non-dispersive, non-chiral case are unity and $\lambda$ (see (\ref{dimgoveqn}) with $\gamma=0$). In the chiral case, where the gyroscopic coupling leads to dispersion, the displacement components are displayed in  Figs. \ref{fig2a} and \ref{fig2b}, and the profiles associated with non-chiral wavefronts are also shown. The characteristic lines showing these wavefronts, $x\pm t= \mbox{const}$  and $x\pm \lambda t= \mbox{const}$, are marked on both Figs. \ref{fig2a} and \ref{fig2b} with dashed and dotted lines respectively. Note that all the wave motion for both components is bounded by the characteristic corresponding to the higher wave speed.

\subsection{A temporal chiral interface of high gyricity}
\label{tempch}

We assume that at time $t=0$, the gyricity switches from $0$ to $\Ga \gg 1.$ At time $t=d >0$, the gyricity switched back to zero. In this case, we consider the temporal interface $0 < t < d$. In particular, we are interested in the values of $\Ga$ and $ d$ such that $\Ga \gg 1$, while the product $\Ga d$ remains finite.  

It has been demonstrated above that for high values of $\Ga$, the vibrations of the chiral system, discussed here, can be considered on a cross-section, with fixed $x$, and we use a vector function 
$\BY(\Ga, t)$, where the $x-$dependence is omitted.

Consider the following initial value problem
\beq
\fr{\prt }{\prt t} \Big( \BCM (\Ga t) \fr{\prt}{\prt t} \BY(\Ga, t) \Big)  = {\bf 0},
\eequ{ceq20} 
\beq
\BY(\Ga, 0) = \Bf = \begin{pmatrix} f_1 \cr f_2 \end{pmatrix}, ~ \BCM ( 0) \fr{\prt}{\prt t} \BY(\Ga, 0) = \fr{\prt}{\prt t} \BY(\Ga, 0) = \Bg=\begin{pmatrix} g_1 \cr g_2 \end{pmatrix}, 
\eequ{ceq21}
with $\Bf, \Bg$ being the  vectors representing the initial displacements and initial velocities at $t=0.$ 

The solution of \eq{ceq20}, \eq{ceq21} has the form
\beq
\BY(\Ga, t) = \Ga^{-1} \BR (\BI - \BCM^T(\Ga t)) \Bg + \Bf,
\eequ{ceq22}
where $\BI$ is the $2\times 2$ identity  matrix.
We also note that
\beq
\fr{\prt}{\prt t} \BY(\Ga, t) = \BCM^T(\Ga t) \Bg,
\eequ{ceq23}
and hence the momentum remains constant for all positive $t$
\beq
\BCM (\Ga t) \fr{\prt}{\prt t}  \BY (\Ga, t) =  \Bg = \BCM(0) \fr{\prt}{\prt t}  \BY(\Ga, 0). 
\eequ{ceq24}
It also follows from \eq{ceq22} that
\beq
\BY(\Ga, d) - \BY(\Ga, 0)   = \Ga^{-1} \BR (\BI - \BCM^T(\Ga d))  \fr{\prt}{\prt t}  \BY (\Ga, 0),
\eequ{ceq25}
which shows that the displacement $\BY$ has a discontinuity of a small amplitude across the chiral interface.

 Such an interface also provides the rotational coupling between the longitudinal and transverse vibrations.  

In particular, if $\Ga$ and $d$ are chosen in such a way that $\Ga d = (2n - 1) \pi$, where $n$ is positive integer, then \eq{ceq25} becomes
\beq
\Big[   \BY(\Ga, t) \Big]_{t=0}^{t=d} =  2 \Ga^{-1} \BR \Bg,
\eequ{ceq26}
whereas for the case of  $\Ga d = 2n \pi$ we have the ideal temporal interface, across which both the displacement and the momentum are continuous
\beq
\Big[  \BY \Big]_{t = 0}^{t=d} = {\bf 0}, ~~  \Big[ \BCM \fr{\prt}{\prt t}  \BY \Big]_{t = 0}^{t=d} = {\bf 0}.
\eequ{ceq27}
We also note that, with the choice of $\Ga d = \pm \fr{\pi}{2} + 2\pi n $, where $n$ is positive integer, the equation \eq{ceq25} yields
\beq
\Big[   \BY(\Ga, t) \Big]_{t=0}^{t=d} =  \Ga^{-1} (\BR \pm \BI) \Bg.
\eequ{ceq27a}

   
\subsection{An imperfect temporal interface}

Assume that outside the time-interval $t \in ( T^*, T^* + d)$ the gyricity parameter $\Ga$ is zero, and hence the equation  \eq{ceq17} has the solution $\BU(x, t)$, which describes two uncoupled waves propagating with constant speeds $c_1$ (longitudinal) and $c_2$ (transverse), respectively, which are classical D'Alembert's travelling waves. We note that $\BCM (0)= \BI$.

However, it is assumed that within the interval $(  T^*, T^* + d)$ the gyricity parameter $\Ga$ takes a large value, and within that time-interval, on a cross-section with fixed $x$, the solution is described in Section  \ref{tempch}. The continuity of the displacement and of the momentum (see \eq{ceq19}) are set at $t=T_*$ and $t=T^*+d$. 

For small $d$ and large $\Ga$, we assume that $\Ga d = (2n-1) \pi$, where $n$ is positive integer. 
Then it follows from \eq{ceq24}--\eq{ceq26} that across the thin temporal chiral interface the field $\BU$ is discontinuous and the following interface conditions hold:
 
\beq
\Big[   \fr{\prt}{\prt t} \BU(x, t)  \Big]_{t=T_*}^{t=T^*+d} = 0, ~ \Big[   \BU(x, t) \Big]_{t=T_*}^{t=T_*+d} =  2 \Ga^{-1} \BR \fr{\prt}{\prt t} \BU(x, T_*)  .
\eequ{ceq28}

The above interface conditions represent the imperfect temporal interface, across which the displacement vector has the discontinuity. For large $\Ga$ and a finite magnitude of the velocity vector, the right-hand side in the second interface condition \eq{ceq28} is small. On the other hand, when the incident wave approaches the interface with the instantaneous velocity, which is large, then the small coefficient $\Ga^{-1}$ may counterbalance it, and hence a finite magnitude in the displacement jump may be observed across the imperfect temporal interface.

\section{Examples}
\label{S_examples}

Here, we consider three examples for a system of temporal chiral interfaces at $t= T, 2 T, \dots,n T, \dots$, with integer $n$, combined with the Cauchy problems for two wave equations  \eq{ceq1}, where $0 < c_2 < c_1 =1$, and the chirality parameter is large, i.e. $\Ga \gg 1$. In this case, the duration $d$ for each temporal interface is infinitesimally small. 

The initial velocities, at time $t=0$, for both longitudinal and transverse orientations, are assumed to be zero:
\beq
\fr{\prt }{\prt t} u(x, 0) = 0, ~ \fr{\prt }{\prt t} v(x, 0) = 0, 
\eequ{ex_IC1}
while the initial values of the functions $u, v$ may be defined as required.  

Here, we consider three cases:
\begin{enumerate}
\item In the first case, we assume that, in addition to \eq{ex_IC1}, the initial conditions for \eq{ceq1} are chosen in such a way, that 
\beq
u(x, 0) = 0, ~~ v(x,0) = \Phi(x),  
\eequ{ex_IC2}
where $\Phi(x)$ is a smooth even function, exponentially vanishing at infinity. In this case, in the time interval $0< t < T$ only the transverse component $v$ of the displacement is non-zero, whereas $u \equiv 0$. 

\item In the second case,  in addition to \eq{ex_IC1}, we choose the initial conditions for \eq{ceq1} in such a way that 
\beq
u(x, 0) = \Psi(x), ~~ v(x,0) = 0,  
\eequ{ex_IC3}
where $\Psi(x)$ is a smooth even function, exponentially vanishing at infinity. Here, in the time interval $0< t < T$ only the longitudinal component $u$ of the displacement is non-zero, whereas $v \equiv 0$.

\item  In the third case,  both $u$ and $v$ are non-zero i.e.
\beq
u(x, 0) = \Psi(x), ~~ v(x,0) = \Phi(x).  
\eequ{ex_IC4}

\end{enumerate}

At every temporal interface $t = T n$, there is a coupling governed by the transmission conditions \eq{ceq28}. This gives an additional feature in the dynamic response of the elastic system, which incorporates a chiral temporal interface, that couples longitudinal and transverse displacements.  

\subsection{The case of the dominant transverse displacement}

For the case when the initial conditions are chosen to have the form \eq{ex_IC1}, \eq{ex_IC2}, and $|\Ga|$ is a large parameter, it can be shown that the transverse displacement $v$ becomes dominant. The coupling, which occurs at chiral temporal interfaces, yields the solution, which includes a new pattern that incorporates a split of the wave front and the longitudinal displacement, as demonstrated in the closed form analytical representation below.

\begin{figure}
\centering
\includegraphics[height=9.0cm]{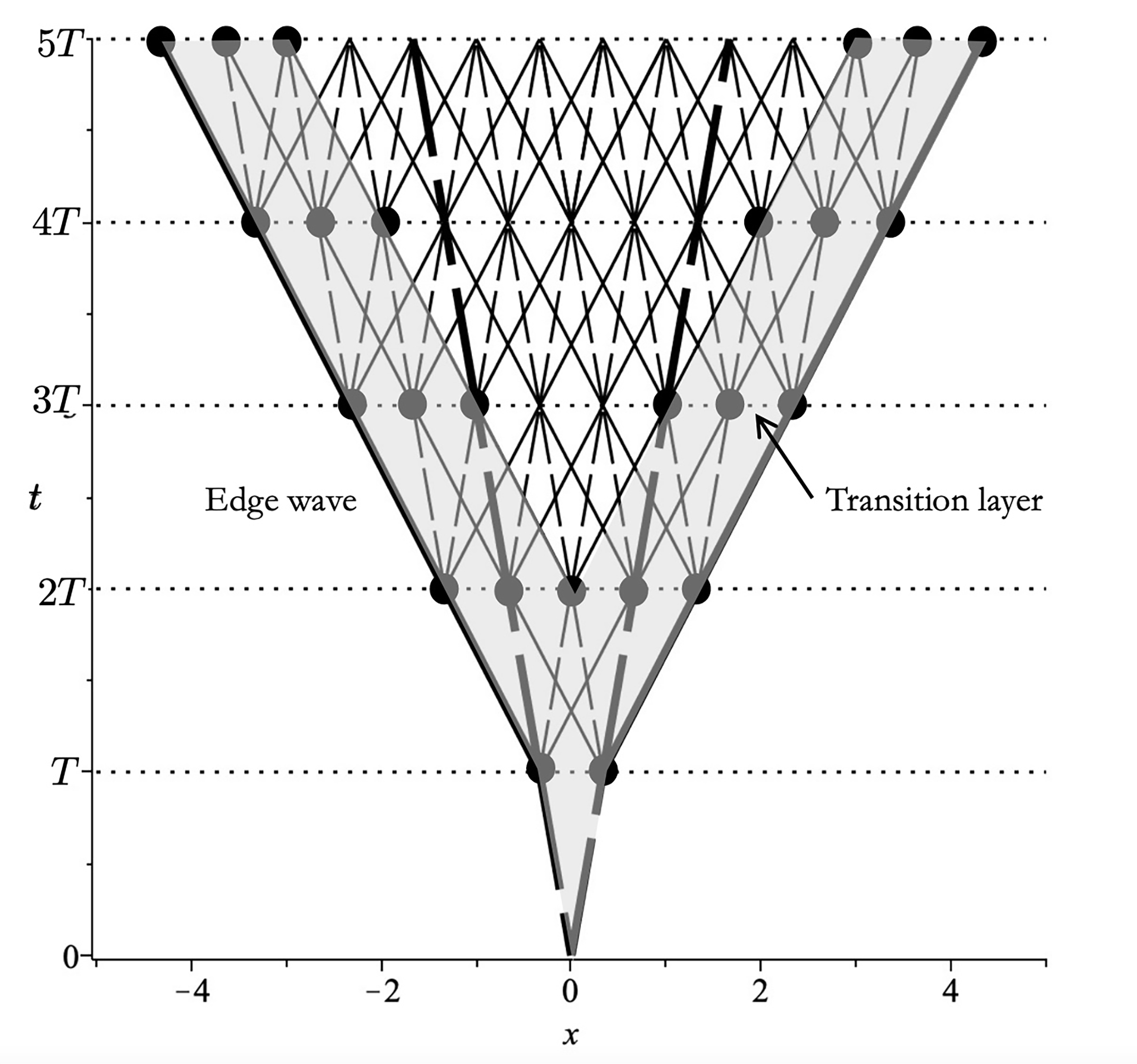}
\caption{\footnotesize Wave split at temporal chiral interfaces. Two families of characteristics are shown, corresponding to wave speeds $c_1 =1$ (solid lines) and $ c_2 = 1/3$ (dashed lines).  The edge wave is shown and also the 'transition layer' where an incomplete number (less than four) of characteristics meet an interface from below.} 
\label{figchar}
\end{figure}

Using D'Alembert's representation for the solution of the Cauchy problem on each of the temporal intervals $(n-1) T < t < n T$, for positive integer $n$, we deduce
\beq
u(x, t) \equiv 0, ~ v(x, t) = \fr{1}{2} \Big( \Phi(x + c_2 t) + \Phi(x-c_2 t)  \Big) ~~\mbox{when} ~0 < t < T,
\eequ{C1_1}
\begin{eqnarray}
u(x,t) &=& \fr{c_2}{2 \Ga} \Big( \Phi'(x+c_1 t - T(c_1 - c_2)) + \Phi'(x-c_1 t + T (c_1 + c_2)) \nonumber \\
&& - \Phi'(x+c_1 t - T(c_1 + c_2)) - \Phi'(x-c_1 t + T(c_1-c_2)) \Big), \nonumber \\
v(x, t) &=&  \fr{1}{2} \Big( \Phi(x + c_2 t) + \Phi(x-c_2 t)  \Big) ~~\mbox{when} ~T < t < 2T, \label{C1_2}
\end{eqnarray}
\begin{eqnarray}
u(x, t) &=& \fr{c_2}{2 \Ga} \Big(   \Phi'(x+c_1 t - T(c_1 - c_2)) + \Phi'(x-c_1 t + T(c_1 + c_2)) \nonumber \\
&& - \Phi'(x+c_1 t - T(c_1 + c_2)) - \Phi'(x-c_1 t + T(c_1 - c_2)) \nonumber \\
&& + \Phi'(x+c_1 t - 2 T(c_1 - c_2)) + \Phi'(x-c_1 t + 2 T (c_1 + c_2)) \nonumber \\
&& - \Phi'(x+c_1 t - 2 T(c_1 + c_2)) - \Phi'(x-c_1 t + 2 T (c_1 - c_2))   
\Big) ,\nonumber \\
v(x, t) &=& \fr{1}{2} \Big( \Phi(x+c_2 t) + \Phi(x-c_2 t) \Big) \nonumber \\
&& - \fr{c_1 c_2}{2 \Ga^2} \Big( \Phi''(x+c_2 t + T(c_1 - c_2)) + \Phi''(x-c_2 t + T(c_1 + 3c_2)) \nonumber \\
&&- \Phi''(x+c_2 t - T(c_1 + c_2)) - \Phi''(x-c_2 t - T(c_1 - 3 c_2)) \nonumber \\
&&- \Phi''(x+c_2 t + T(c_1 - 3 c_2)) - \Phi''(x-c_2 t + T(c_1+c_2)) \nonumber \\
&&+ \Phi''(x+c_2 t - T(c_1 + 3 c_2)) + \Phi''(x-c_2 t - T(c_1 - c_2)) \Big) ~~ \nonumber \\
&&  \mbox{when} ~2T < t < 3T,  \label{C1_3}
\end{eqnarray}
which can be continued further, using the standard D'Alembert representation on each of the temporal interfaces, to show that the transverse displacement  $ v(x, t) = \fr{1}{2} \Big( \Phi(x + c_2 t) + \Phi(x-c_2 t)  \Big) $ remains dominant at all times, with the magnitude of waves associated with the longitudinal displacement $u(x, t)$ being of order $O(|\Ga|^{-1})$.

\begin{figure}[H]
	\centering
	\begin{subfigure}{.4\textwidth}
		\includegraphics[width=\textwidth]{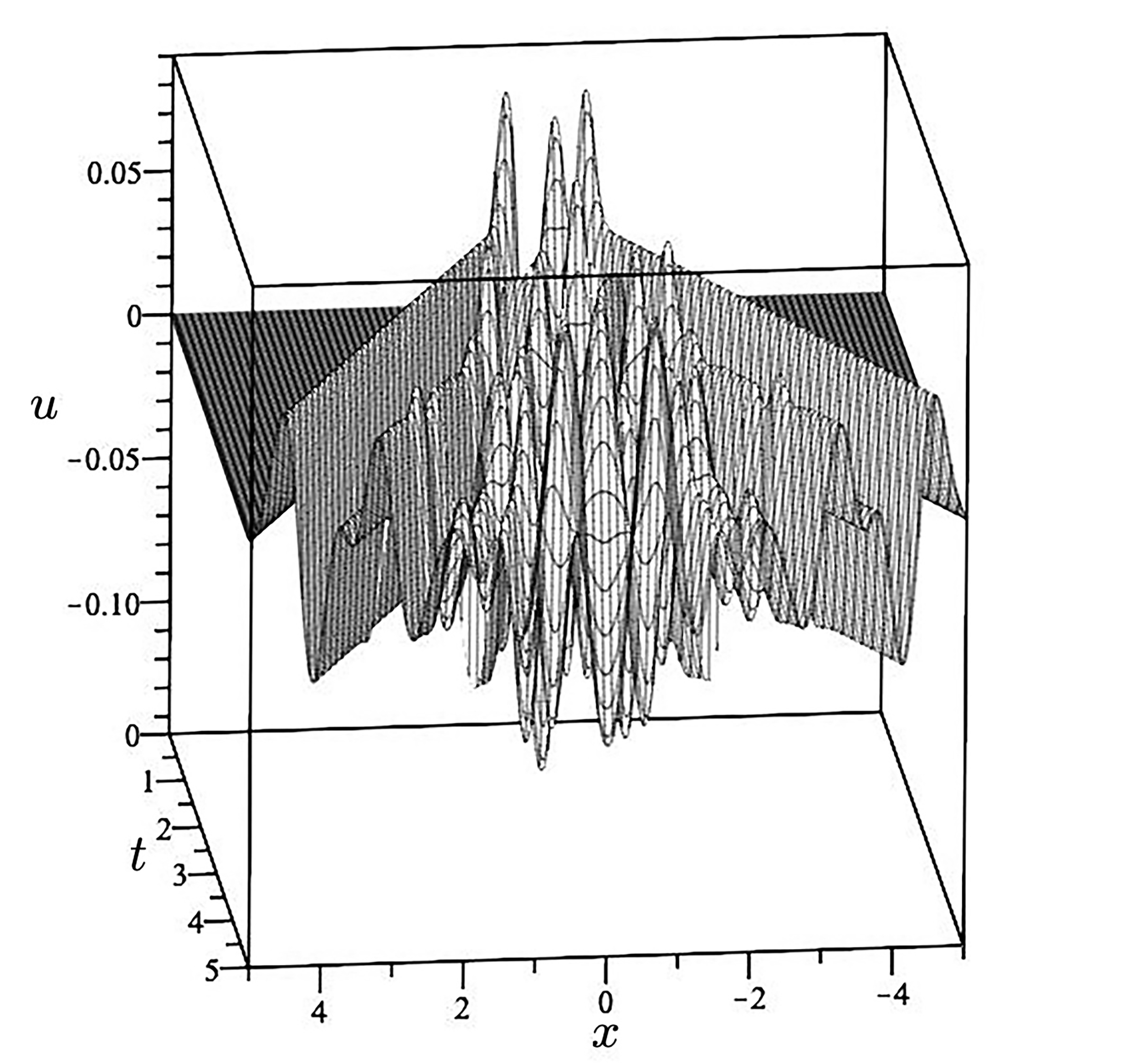}
		\caption{\footnotesize $u(x,t)$}
	\end{subfigure}
	~~~\begin{subfigure}{.4\textwidth}
		\includegraphics[width=\textwidth]{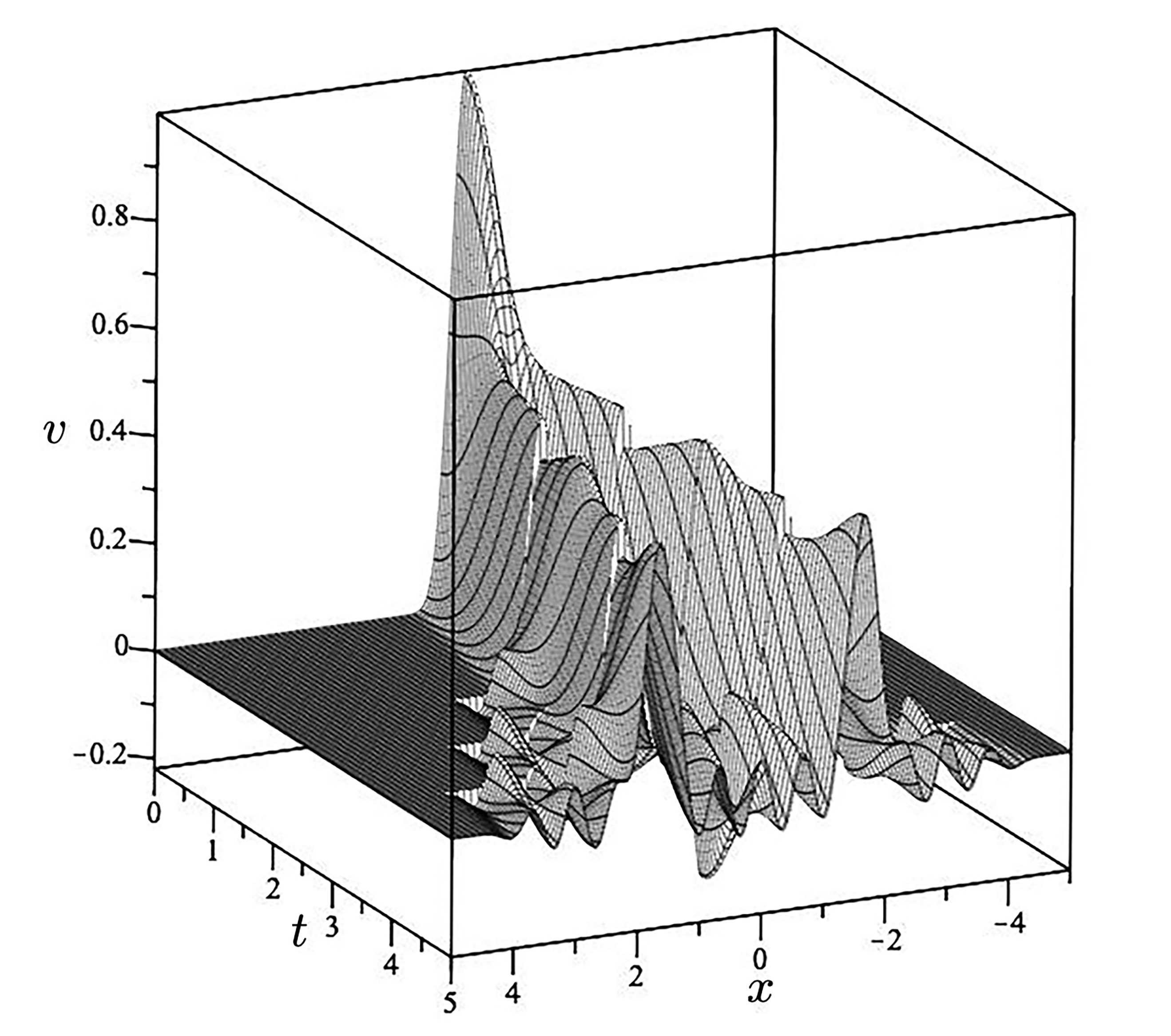}
		\caption{\footnotesize $v(x,t)$}
	\end{subfigure}
	\\
	\begin{subfigure}{.4\textwidth}
		\includegraphics[width=\textwidth]{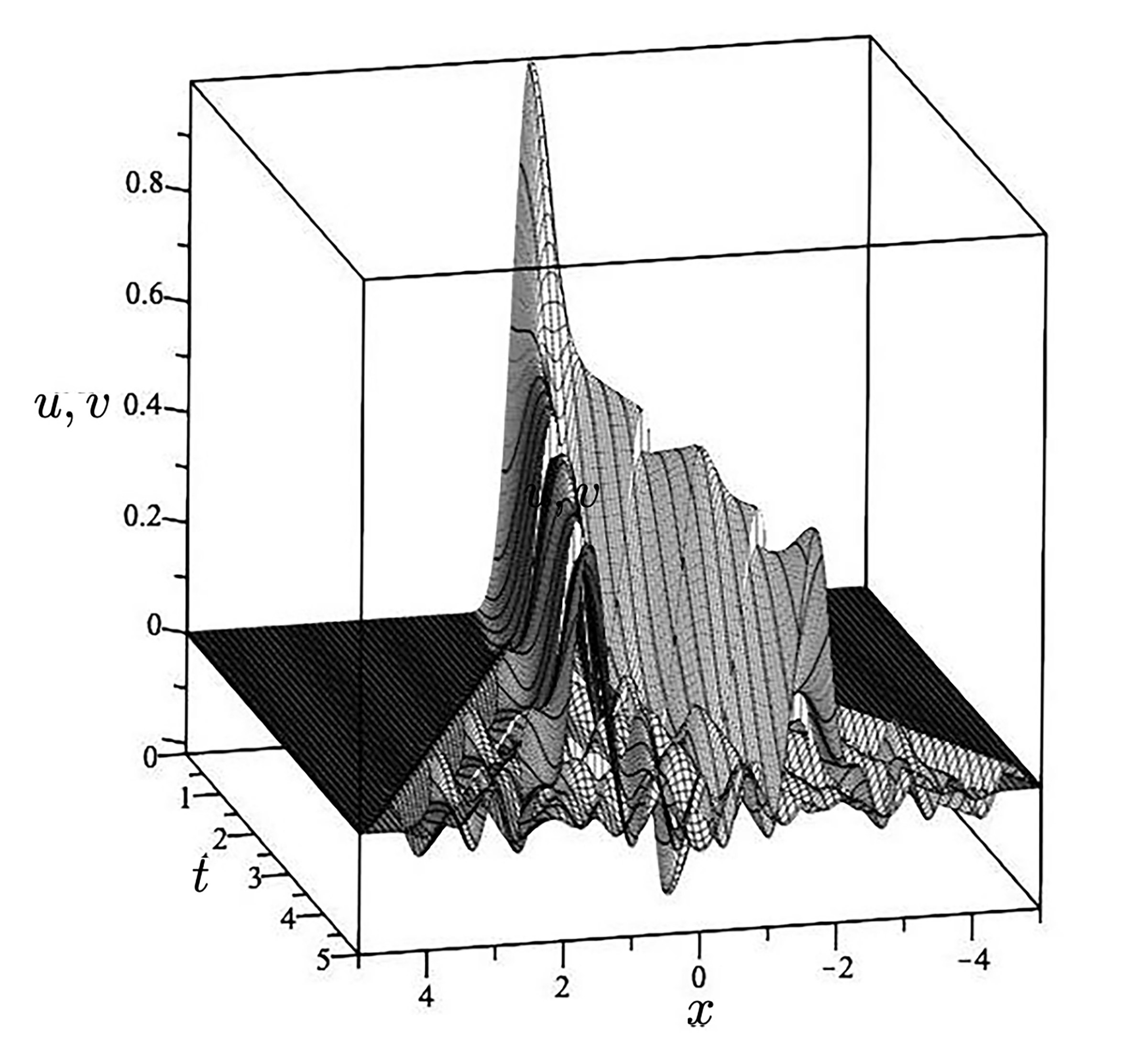}
		\caption{\footnotesize two surfaces, corresponding to $u(x,t)$ and $v(x,t)$, are plotted on the same diagram }
	\end{subfigure}
~~~\begin{subfigure}{.37\textwidth}
		\includegraphics[width=\textwidth]{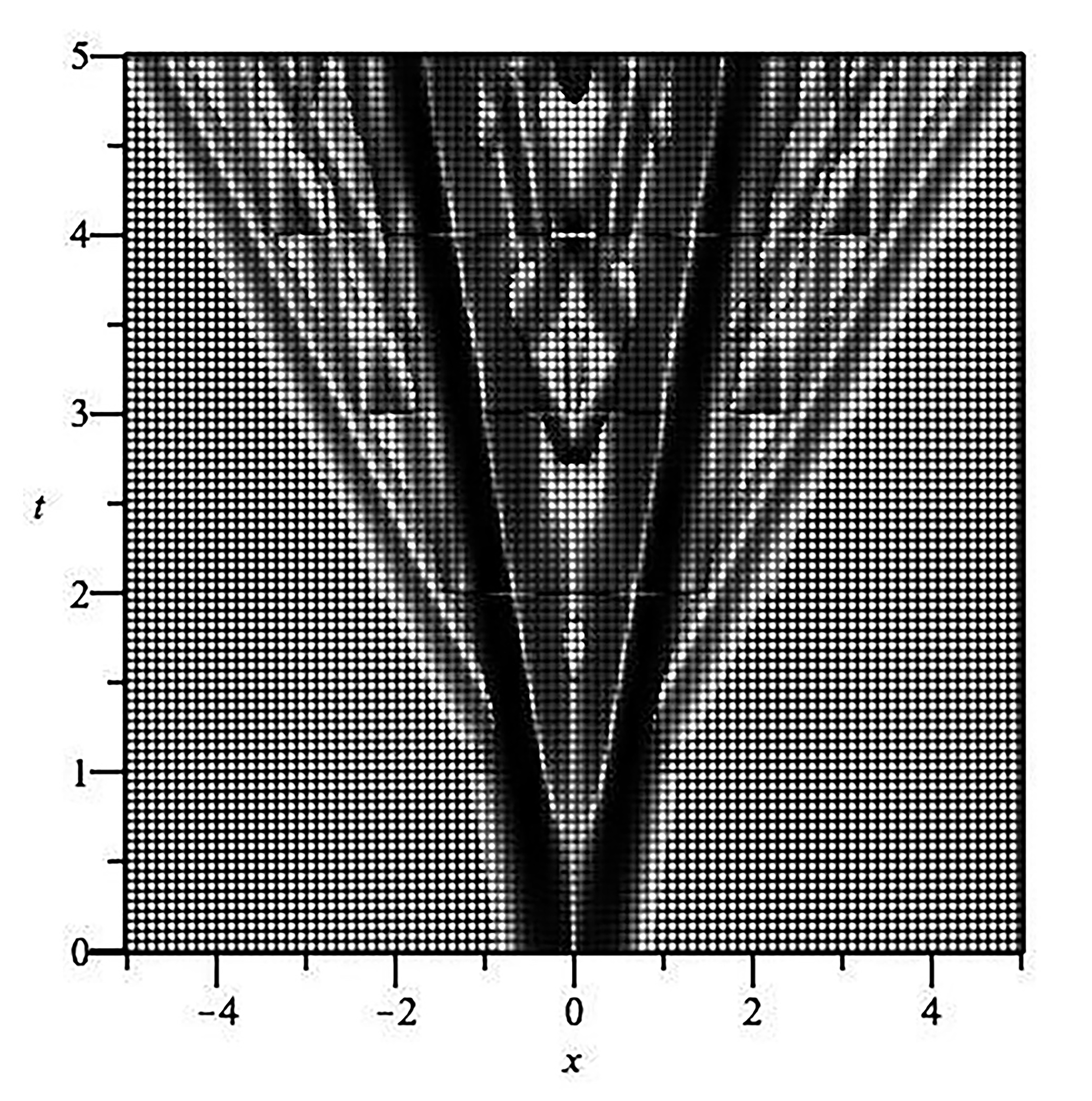}
		\caption{\footnotesize the same surfaces, $u(x,t)$ and $v(x,t)$, as in part (c) - view from above}
	\end{subfigure}	
	\caption{\footnotesize The solution of the Cauchy problem for the case 1 (dominant transverse displacement).
The following values of the parameters were used in the computations: $c_1 = 1, c_2 = 1/3, \Ga = 10, A=10.$}
\label{Case1_soln}
\end{figure}

Numerical computations are presented here for the case when $T=1, c_1=1, c_2 = 1/3$. Also, the right-hand side in the initial condition \eq{ex_IC2} is given as $\Phi(x) = \exp(-A x^2),$ with $A$ being a positive constant. 
The diagram of characteristics  is shown in Fig. \ref{figchar} (compare with Fig. \ref{Fig2}).

In this case, the wave propagates with the speed $c_2$ during the time interval $0 < t < T,$ and the equation of the characteristics is
$t = |x|/c_2.$ 

At time $t=T$ the coupling, induced by the chiral interface, will lead to the formation of two families of waves, propagating with speeds $c_1=1$ and 
$c_2 = 1/3.$
 Taking into account the wave split at the temporal interface, we can write the equations of characteristics at the interval $T < t < 2 T$ as 
\beq
t = |\fr{1}{c_2} |x| - T| +T, ~~t = \fr{1}{c_1} |x \pm T c_2 | + T.
\eequ{eq29}
This process can be repeated further and the resulting diagram of charactreristics is shown in Fig. \ref{figchar}. We note that this diagram does not show the magnitude of the wave. In particular, along some characteristics the wave magnitude may be zero. This is illustrated in the analytical representation of the solution \eq{C1_1}--\eq{C1_3}, as well as in the surface plot of the solution shown in Fig. \ref{Case1_soln}.

\begin{figure}
\centering
\includegraphics[height=9.0cm]{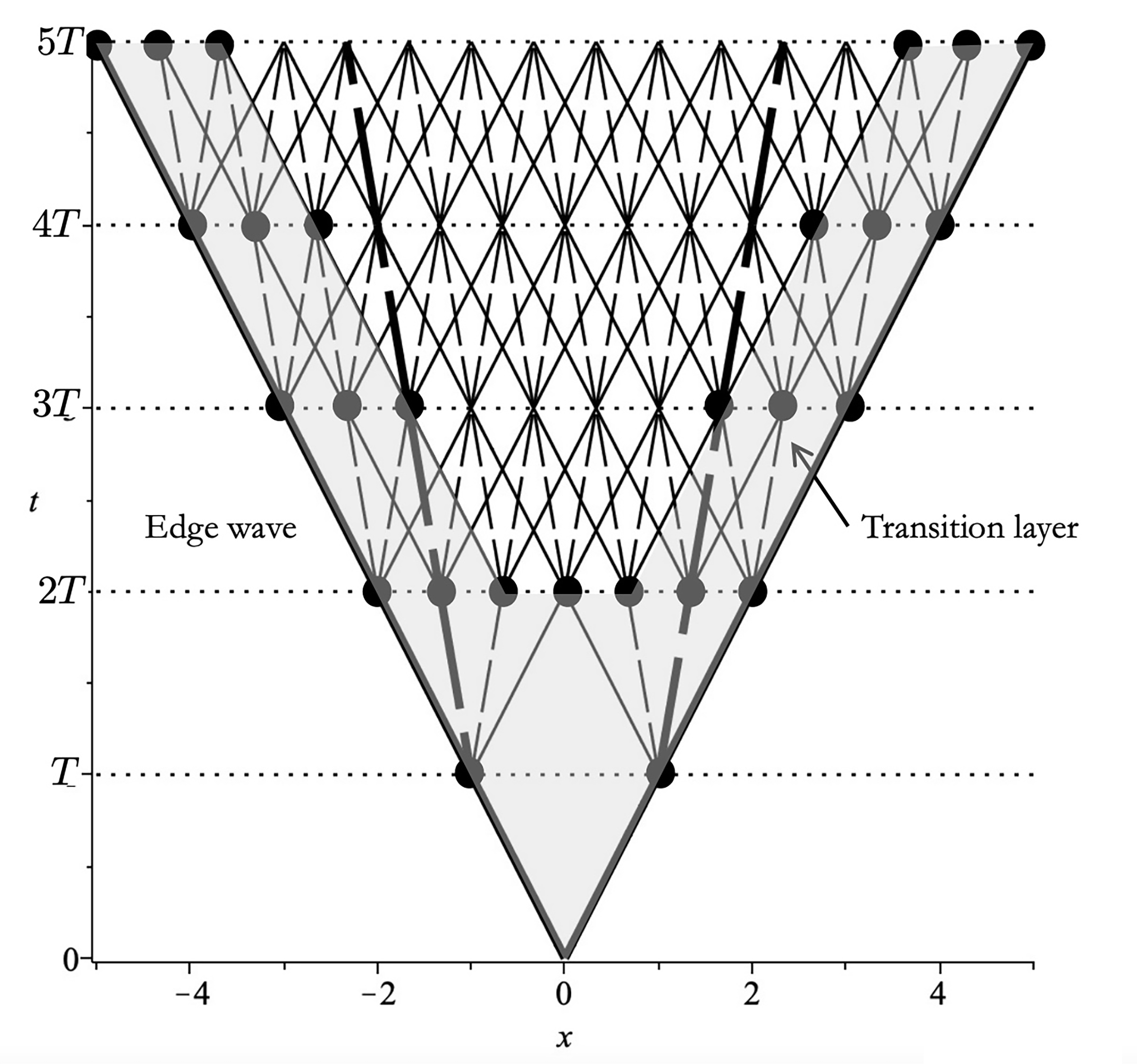}
\caption{\footnotesize The case, when the longitudinal displacement is dominant. Two families of characteristics are shown, corresponding to wave speeds $c_1 =1$ (solid lines) and $ c_2 = 1/3$ (dashed lines).  The edge wave is shown and also the 'transition layer' where an incomplete number (less than four) of characteristics meet an interface from below.} 
\label{figchar1}
\end{figure} 
\subsection{The case of the dominant longitudinal displacement}

\begin{figure}[H]
	\centering
	\begin{subfigure}{.4\textwidth}
		\includegraphics[width=\textwidth]{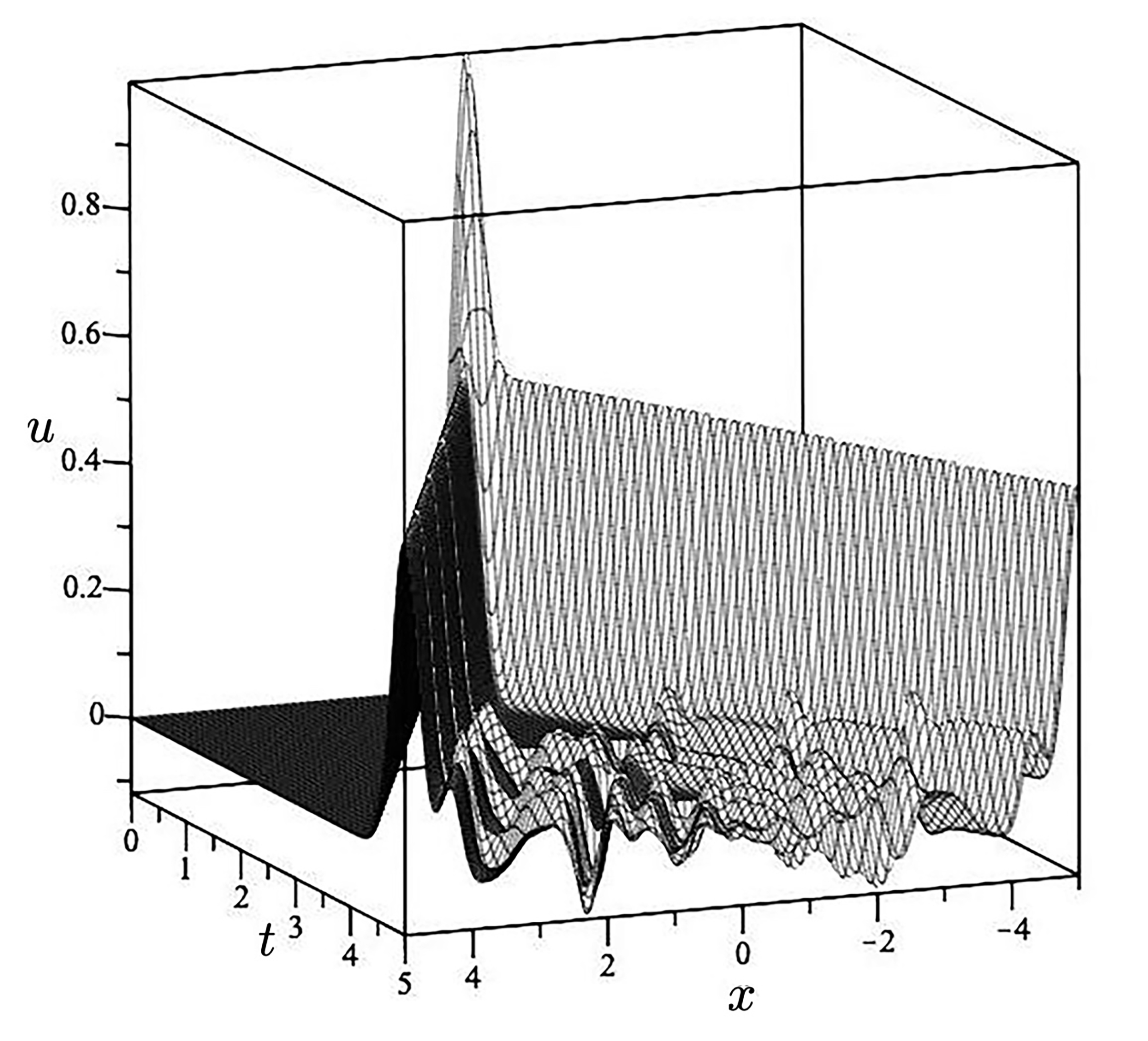}
		\caption{\footnotesize $u(x,t)$}
	\end{subfigure}
~~~	\begin{subfigure}{.4\textwidth}
		\includegraphics[width=\textwidth]{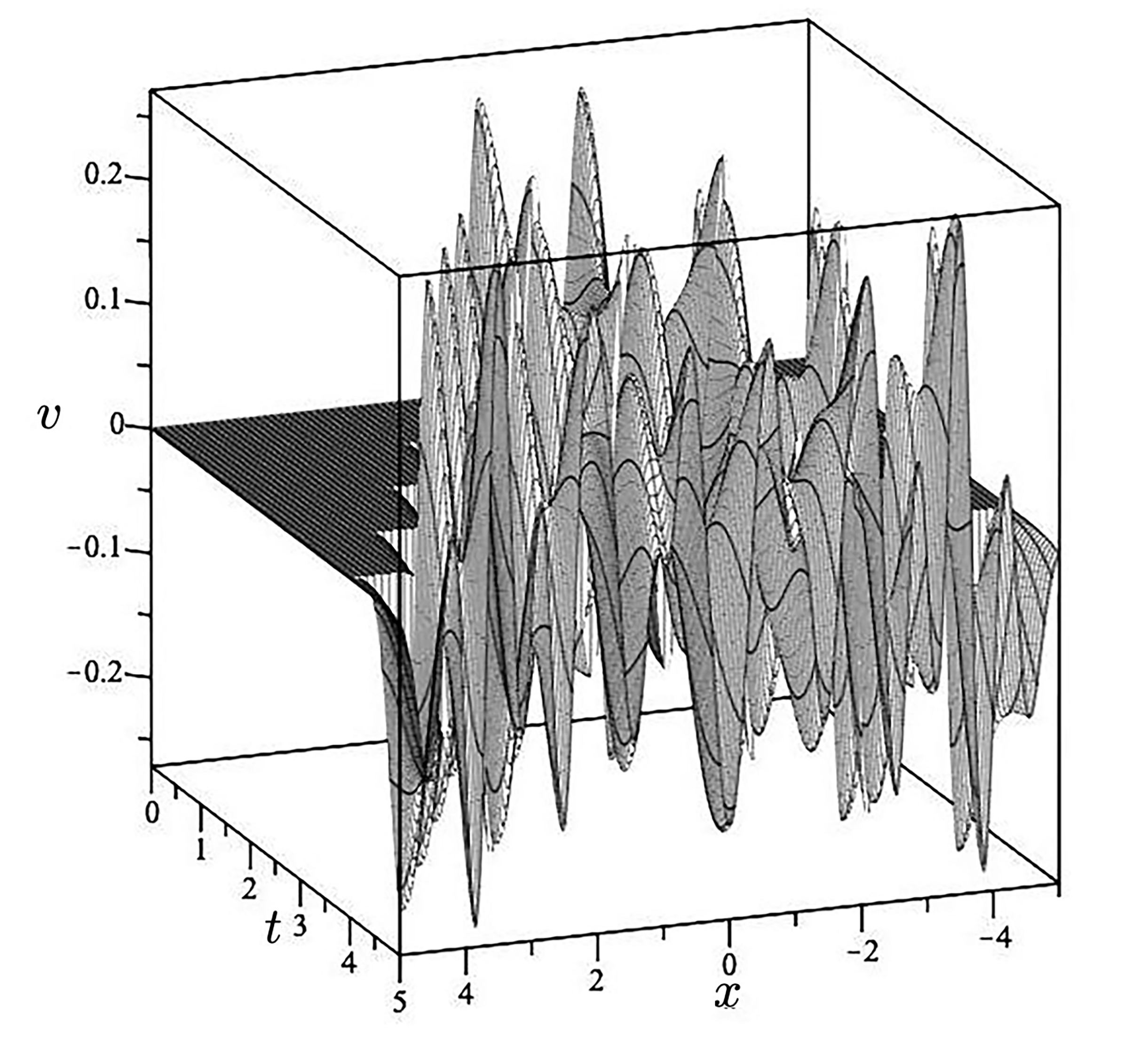}
		\caption{\footnotesize $v(x,t)$}
	\end{subfigure}
	\\
	\begin{subfigure}{.4\textwidth}
		\includegraphics[width=\textwidth]{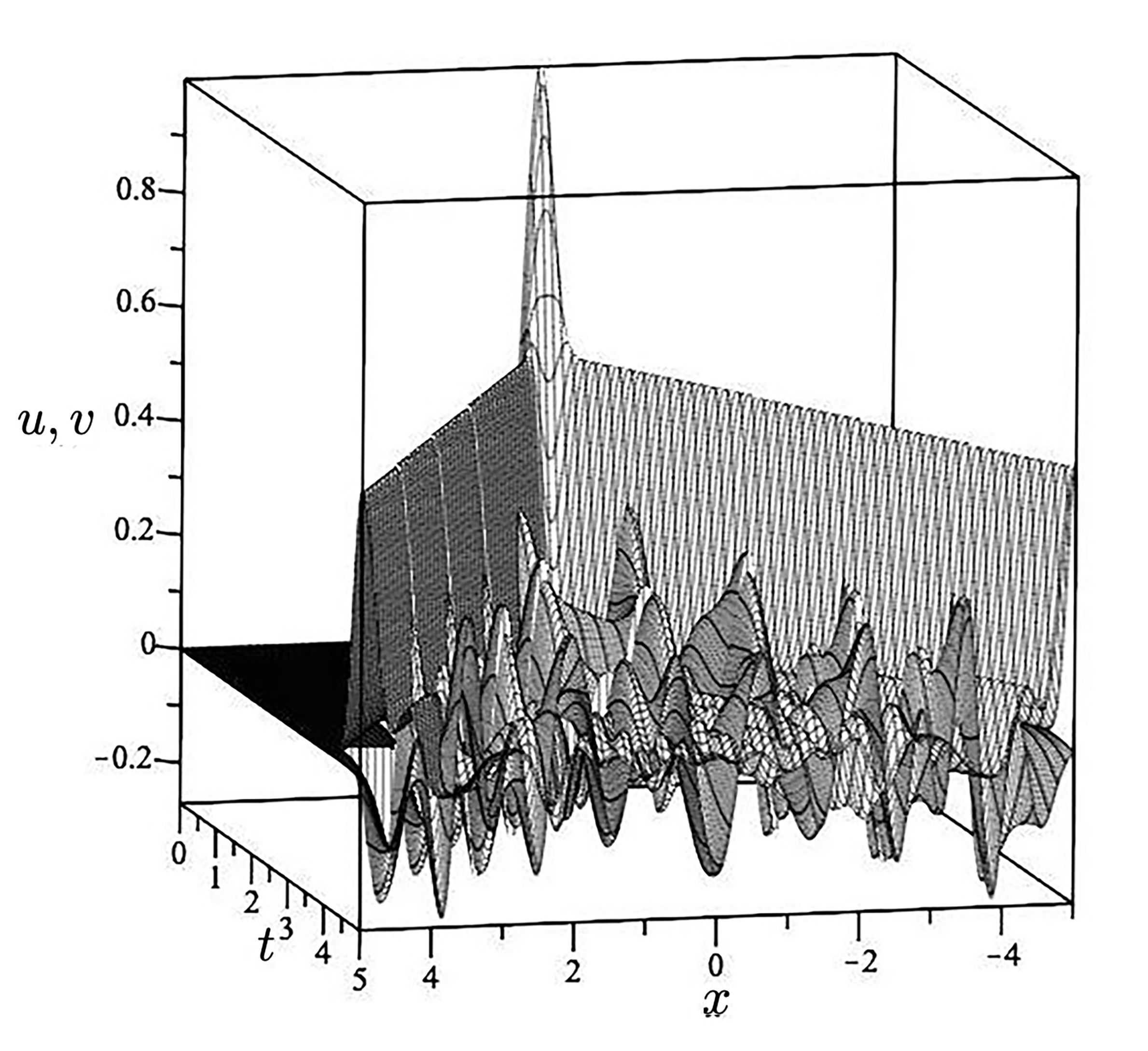}
		\caption{\footnotesize two surfaces, corresponding to $u(x,t)$ and $v(x,t)$, are plotted on the same diagram }
	\end{subfigure}
~~~\begin{subfigure}{.4\textwidth}
		\includegraphics[width=\textwidth]{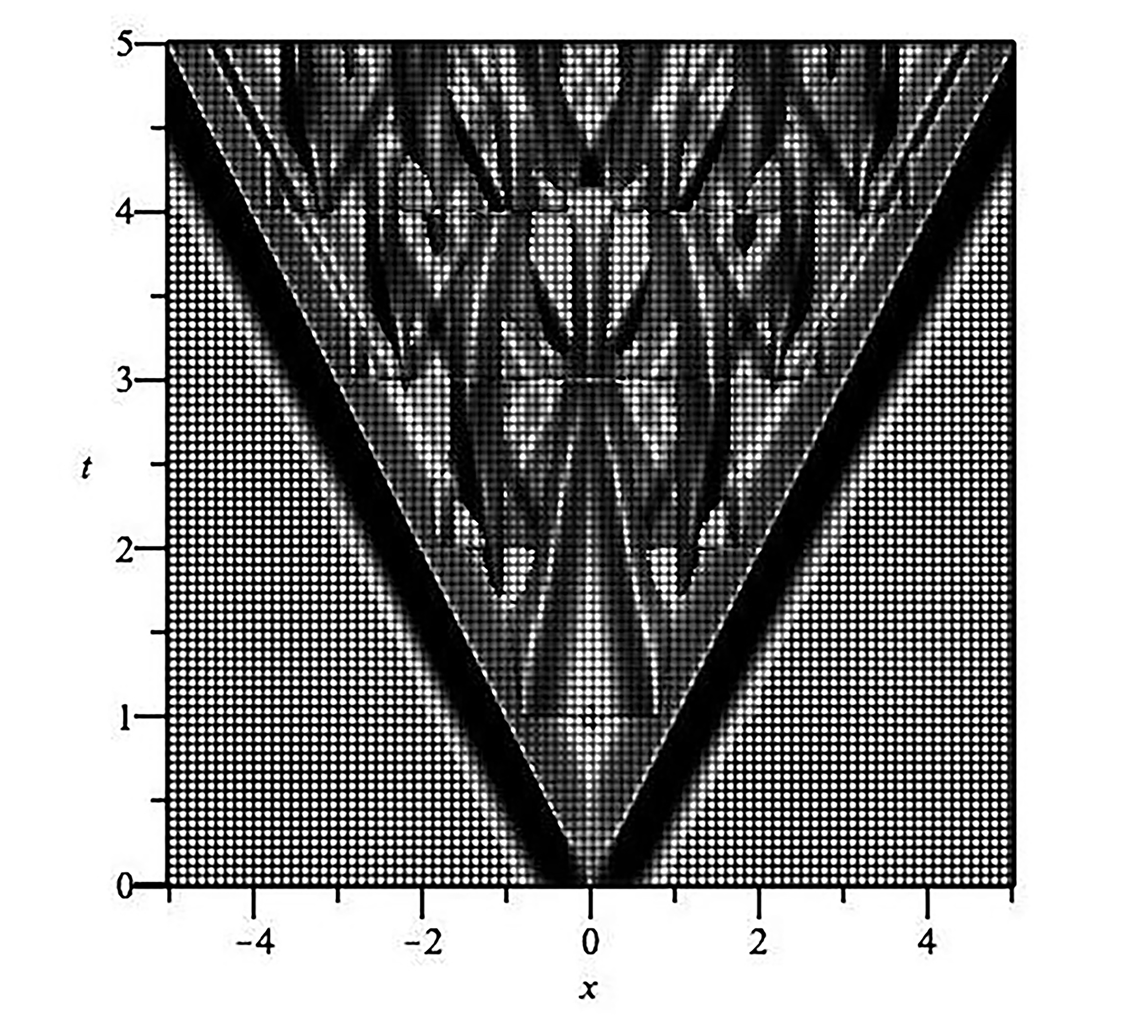}
		\caption{\footnotesize the same surfaces, $u(x,t)$ and $v(x,t)$, as in part (c) - view from above}
	\end{subfigure}	
	\caption{\footnotesize The solution of the Cauchy problem for case 2 (dominant longitudinal displacement).
The following values of the parameters were used in the computations: $c_1 = 1, c_2 = 1/3, \Ga = 10, A=10.$}
\label{Case2_soln}
\end{figure}

Here, we show the example related to the choice of \eq{ex_IC1} and \eq{ex_IC3} as the initial conditions. The right-hand side in the first initial condition \eq{ex_IC3} is given as $\Psi(x) =  \exp(-A x^2),$ with $A$ being a positive constant. As in the previous section, we also choose $T=1, c_1 =1, c_2 = 1/3.$
The diagram of characteristics for this case is shown in Fig. \ref{figchar1}, and the surface plot representing the solution of the Cauchy problem is shown in Fig. \ref{Case2_soln}. We note that the longitudinal displacement is dominant and  the coupling, which occurs at chiral temporal interfaces, yields small wave ripples associated with the transverse displacement.

\subsection{The combined case where both longitudinal and transverse displacements are present at the initial time}

\begin{figure}
\centering
\includegraphics[height=9.0cm]{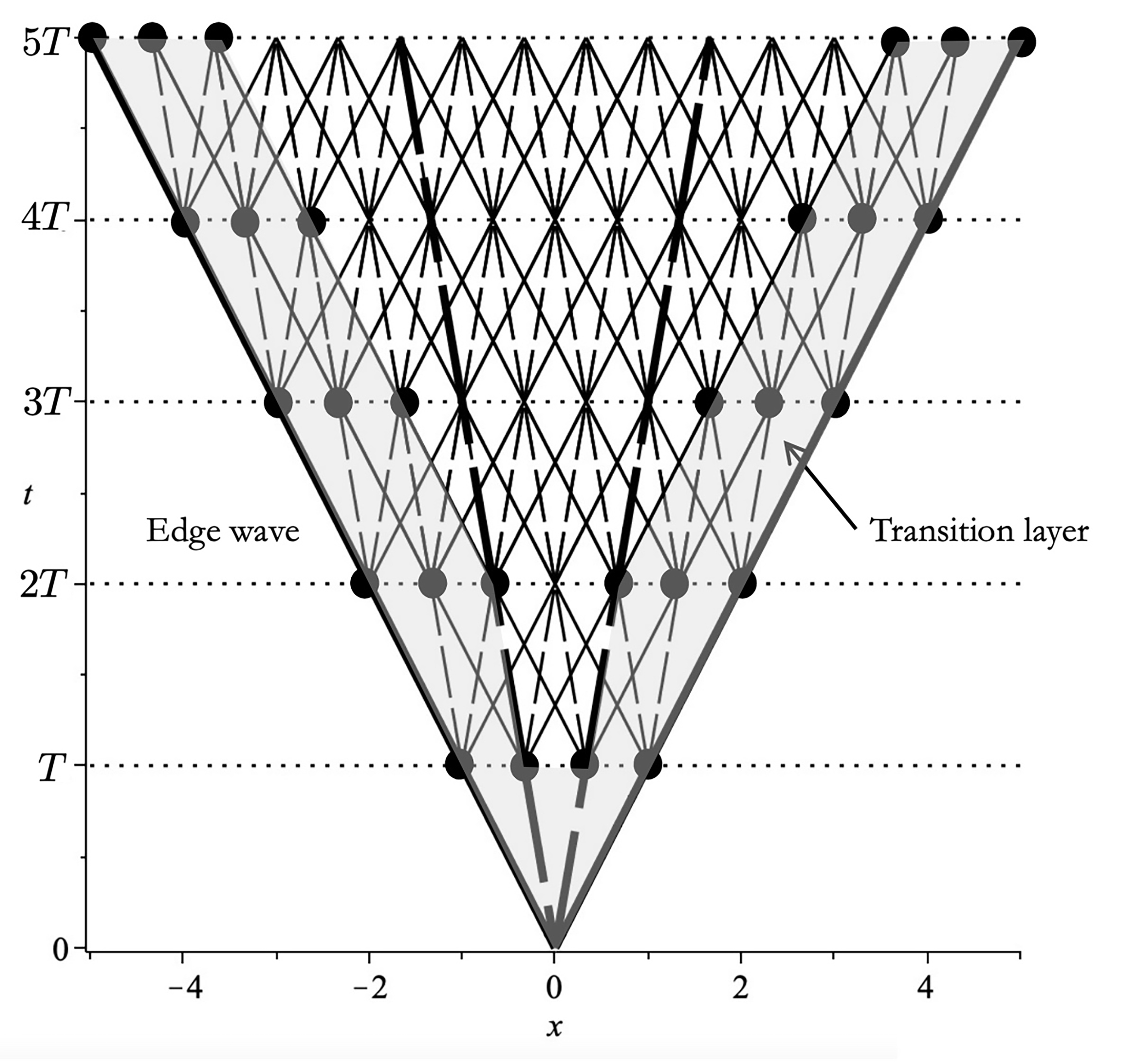}
\caption{\footnotesize The combined case when both the longitudinal and transverse displacements are present in the initial conditions. Two families of characteristics are shown, corresponding to wave speeds $c_1 =1$ (solid lines) and $ c_2 = 1/3$ (dashed lines).  The edge wave is shown and also the 'transition layer' where an incomplete number (less than four) of characteristics meet an interface from below.} 
\label{figchar2}
\end{figure}

Here, the example where both displacements are non-zero at the initial time is illustrated i.e. we choose \eq{ex_IC1}, \eq{ex_IC4} as the initial conditions. The right-hand sides in the initial conditions \eq{ex_IC4} are chosen to be $\Phi(x) = \Psi(x) =  \exp(-A x^2),$ with $A$ being a positive constant. We also choose $T=1, c_1 =1, c_2 = 1/3.$
The diagram of characteristics for this case is shown in Fig. \ref{figchar2}, and the surface plot representing the solution of the Cauchy problem is shown in Fig. \ref{Case3_soln}. Both displacements $u$ and $v$ are present at all times in this case.
However, due to the presence of the chiral temporal interfaces, the wave pattern emerges, which is consistent with the diagram of characteristics of Fig. \ref{figchar2}.  This wave pattern also includes the  ``{transition layer}'', adjacent to the {\em edge wave}, which will both be discussed in the next section.

\subsection{The edge wave and the transition layer}

\begin{figure}[H]
	\centering
	\begin{subfigure}{.4\textwidth}
		\includegraphics[width=\textwidth]{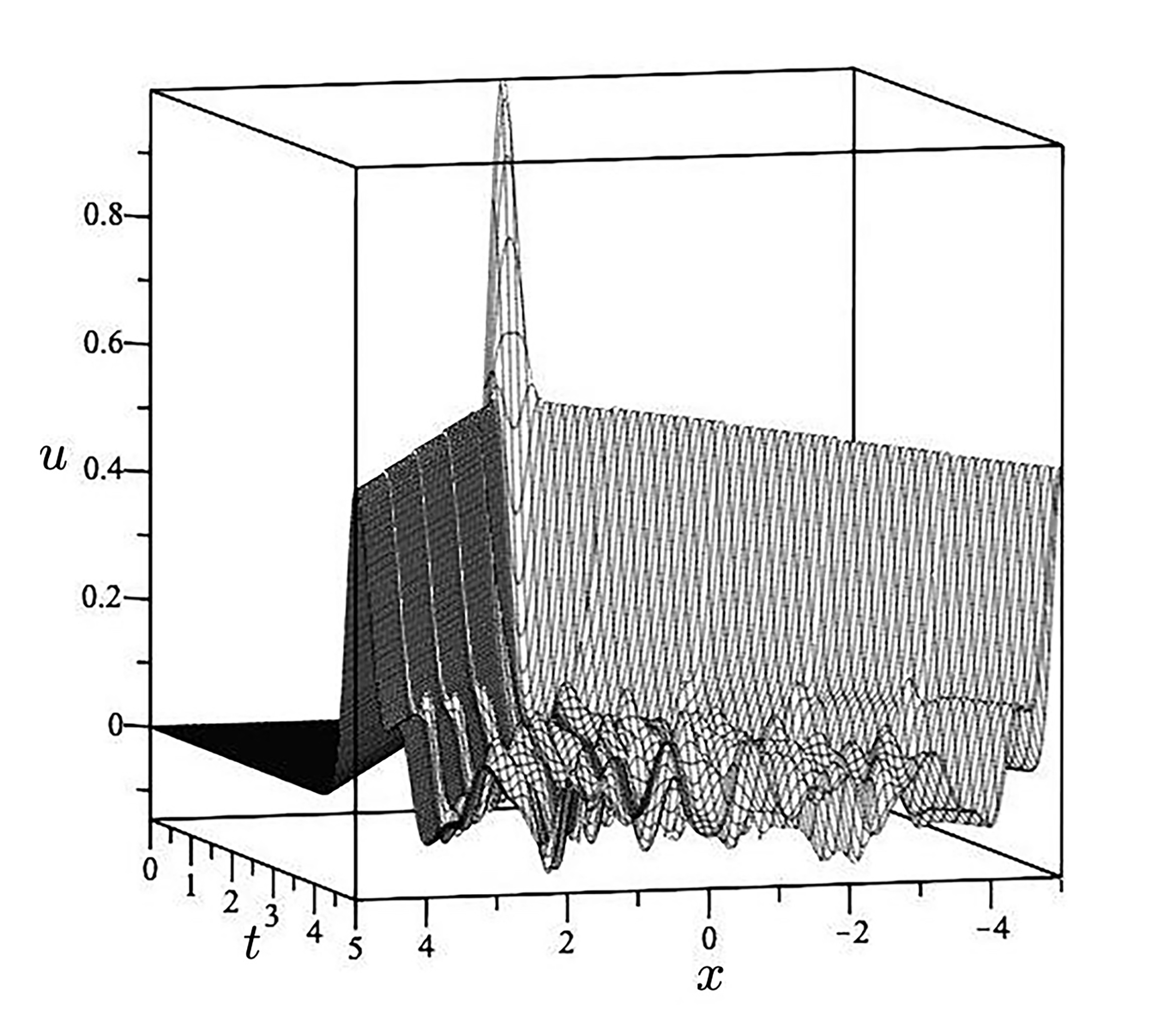}
		\caption{\footnotesize $u(x,t)$}
	\end{subfigure}
~~~	\begin{subfigure}{.4\textwidth}
		\includegraphics[width=\textwidth]{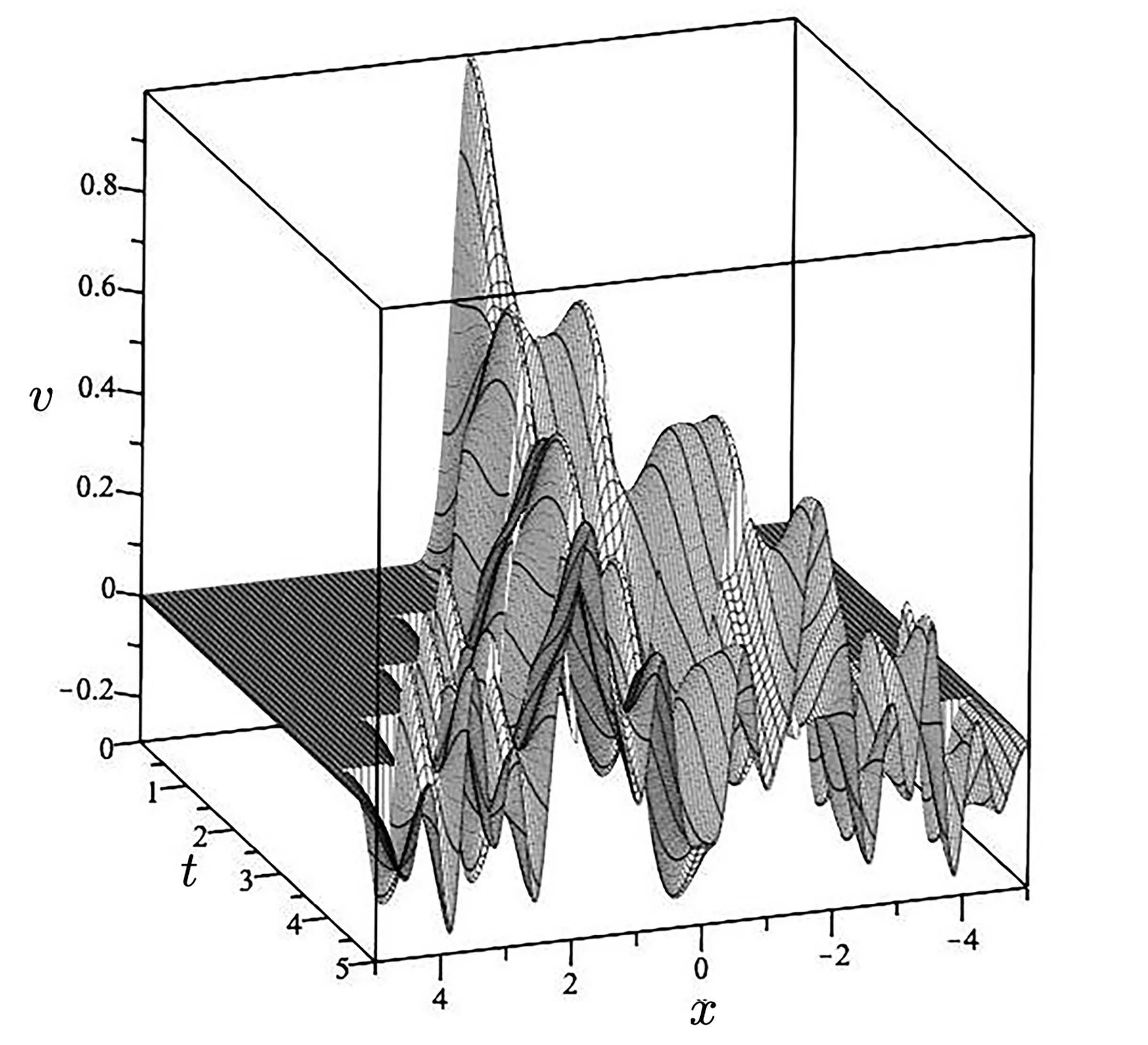}
		\caption{\footnotesize $v(x,t)$}
	\end{subfigure}
	\\
	\begin{subfigure}{.4\textwidth}
		\includegraphics[width=\textwidth]{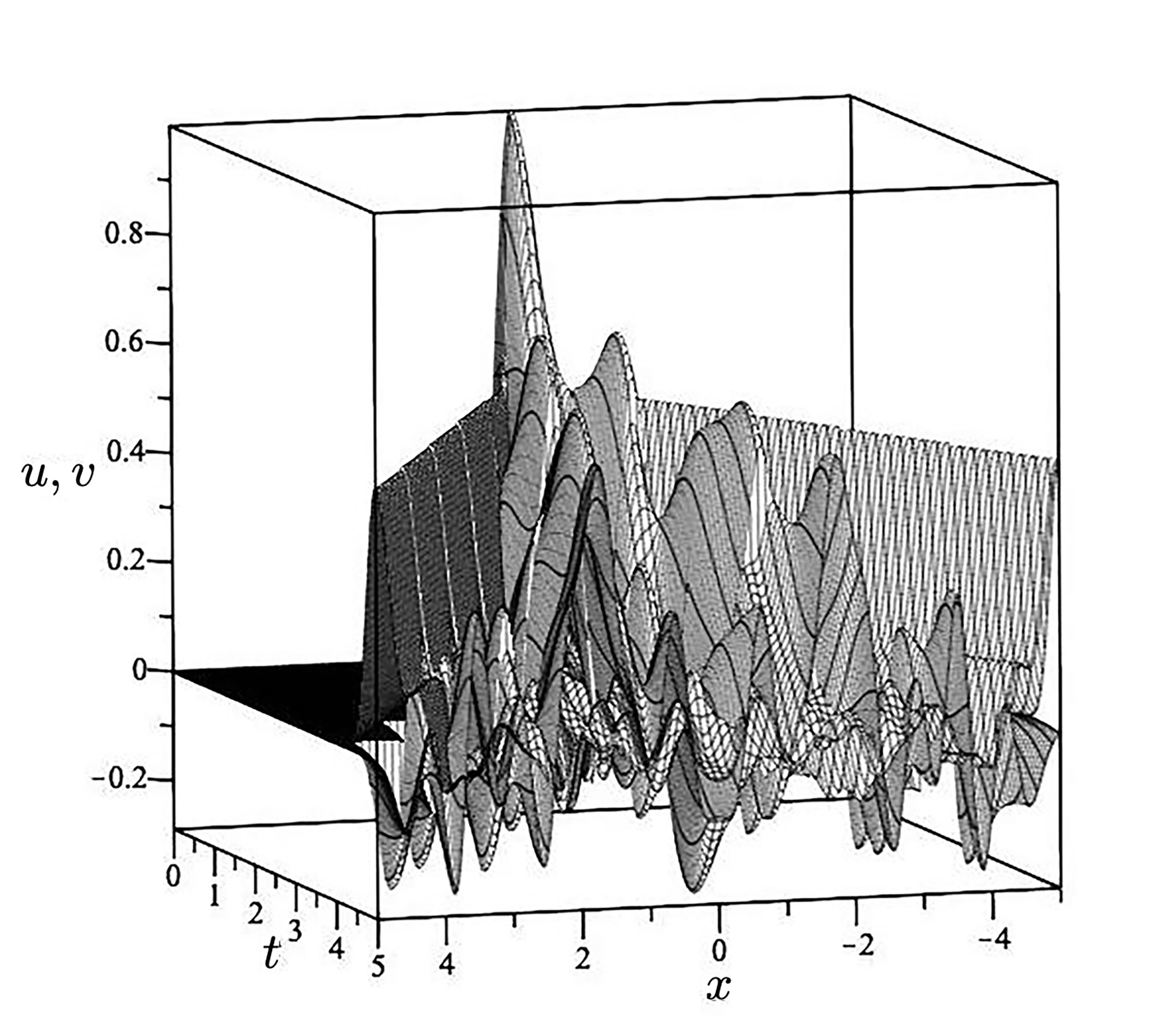}
		\caption{\footnotesize two surfaces, corresponding to $u(x,t)$ and $v(x,t)$, are plotted on the same diagram }
	\end{subfigure}
~~~ \begin{subfigure}{.4\textwidth}
		\includegraphics[width=\textwidth]{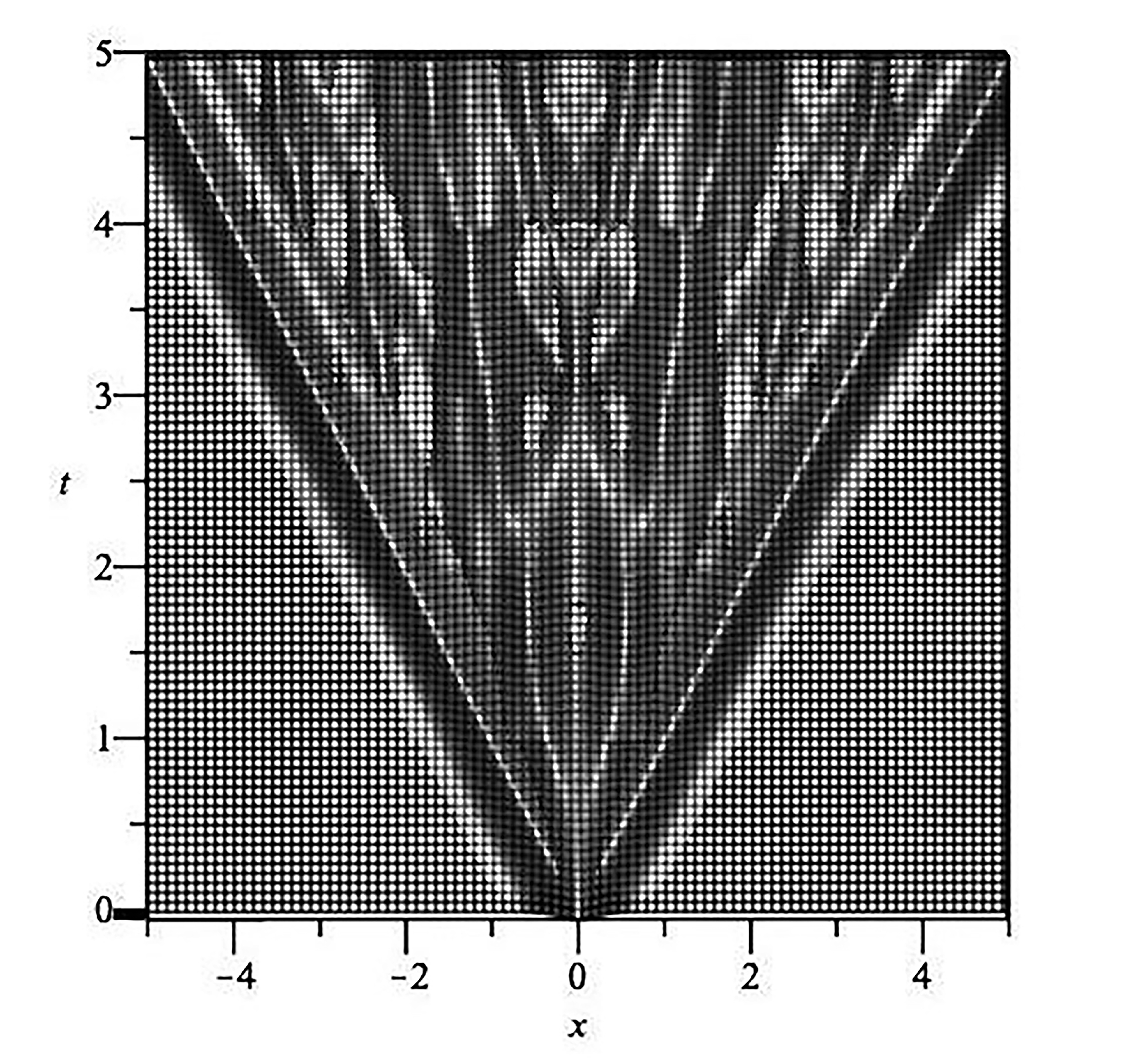}
		\caption{\footnotesize the same surfaces, $u(x,t)$ and $v(x,t)$, as in part (c) - view from above}
	\end{subfigure}	
	\caption{\footnotesize The solution of the Cauchy problem for case 3 (both longitudinal and transverse displacements are present at all times).
The following values of the parameters were used in the computations: $c_1 = 1, c_2 = 1/3, \Ga = 10, A=10.$}
\label{Case3_soln}
\end{figure}

Using the notion of the edge wave introduced in Section \ref{edgewave}, it can be seen from the diagrams of characteristics, shown in Figs. \ref{figchar}, \ref{figchar1}, \ref{figchar2}, that the edge wave in all three cases propagates with the wave speed $c_1$. For suffiiciently large $\alpha$, the magnitude of the edge wave in the first case, where the transverse wave is dominant, is of order $O(|\Ga|^{-1})$, whereas in the remaining two cases the magnitude of the edge wave is of order $O(1)$, as illustrated in Figs. \ref{Case1_soln},    \ref{Case2_soln},    \ref{Case3_soln}.

{With the reference to Figs. \ref{figchar},  \ref{figchar1}, \ref{figchar2}, we also identify a region adjacent to the edge wave boundary as the 'transition layer'. In this region, at every temporal interface an incomplete set of characteristics (less than four) intersects the interface from below.} The width of the transition layer depends on the ratio of the wave speeds $c_1/c_2$.

We also note that   if the chirality parameter $\Ga$ and the interface thickness are chosen in such a way that $\Ga d = 2 \pi n$, with $n$ being positive integer, then the ideal contact conditions \eq{ceq27} hold across the interface,.
In this case, the standard D'Alembert solution will be observed and no coupling will occur at temporal interfaces.

\section{Concluding remarks}

As shown in Sections \ref{TIM} and \ref{TS}, periodic time-variation in the mass density may lead to interesting features of solutions of wave problems. 
More importantly, the model of a temporally stratified medium can be considered as a result of the discretisation of the wave equations with time-dependent  mass density and stiffness coefficients, as in \eq{eq4}.

The advantage of the discretised formulation in the temporally stratified medium is that in every temporal layer, where the coefficients of the governing equations are time-independent, the problem can be solved and a transmission matrix, as in Section \ref{TS}, can be identified. Furthermore, the process can be carried out in iterations to obtain the solution of the Cauchy problem at any given time.
 Although, the diagrams of characteristics,  as in Figs. \ref{Fig2}, \ref{figchar}, \ref{figchar1}, \ref{figchar2}, do not provide information about the wave magnitude, they readily show possible wave patterns corresponding to selected initial conditions as well as the parameters of temporal stratification.


Sections \ref{temporal}, \ref{edge_wave}, \ref{per_initial} describe a scalar problem of vibrations of an elastic string with time-dependent elastic and inertial properties, with the emphasis on the resonance blow-up as $T \to +\infty.$ The second part of the paper is focused on  a vector problem, where the longitudinal and transverse displacements of the vibrating string are dynamically coupled.

Sections \ref{chiral_interface} and \ref{temporal_chiral} introduce coupled governing equations through chirality, together with the notion of imperfect chiral interfaces. Three examples of Cauchy problems in temporally stratified chiral elastic systems are discussed in detail in Section \ref{S_examples}, and they show the wave patterns, as well the coupling process between the longitudinal and transverse vibrations within the string. A particularly interesting feature, displayed in  Figs. \ref{figchar}, \ref{figchar1}, \ref{figchar2}, is the presence of the {\em transition layer} adjacent to the {\em edge wave} in all three cases. Such transition layers also appear in other field pattern geometries: see, for example, Fig.  16 in \cite{MatMil1}.

\section*{Acknowledgements}

GWM thanks the NSF for support from the US National Science Foundation through grant DMS-2107926. ABM, H-MN, and NVM also gratefully acknowledge support from the
 NSF and the hospitality of the University of Utah where this work was initiated.


\end{document}

%% file: macro.tex
\setlength{\arraycolsep}{2pt}
\newcounter{itemnbr}
\newcommand\eq[1] {(\ref{#1})}
\newcommand\eqp[1] {(\ref{eq:#1})}
\newcommand\defn[1] {\ref{def:#1}}
\newcommand\thm[1] {\ref{thm:#1}}
\newcommand\sect[1] {\ref{sec:#1}}
\newcommand\fig[1] {\ref{fig:#1}}
\newcommand\labeq[1] {\label{eq:#1}}
\newcommand\labdefn[1] {\label{def:#1}}
\newcommand\labthm[1] {\label{thm:#1}}
\newcommand\labsect[1] {\label{sec:#1}}
\newcommand\labfig[1] {\label{fig:#1}}
\newcommand\etal{{\it et al.}}
\newcommand\B {\Bigg}
\newtheorem{definition}{Definition}[section]
\newtheorem{lemma}{Lemma}[section]
\newtheorem{theorem}{Theorem}[section]
\newtheorem{corollary}{Corollary}[section]
\newcommand\proof{{\it Proof:}\quad}
\newcommand\remark{{\it Remark:}\quad}
\newcommand\qed{\ \rule[-0.2ex]{0.3em}{1.5ex}}
\newcommand{\bfm}[1]{\mbox{\boldmath ${#1}$}}
\newcommand{\nonum}{\nonumber \\}
\newcommand{\beqa}{\begin{eqnarray}}
\newcommand{\eeqa}[1]{\label{#1}\end{eqnarray}}
\newcommand{\bequ}{\begin{equation}}
\newcommand{\eequ}[1]{\label{#1}\end{equation}}
\newcommand{\Grad}{\nabla}
\newcommand{\Div}{\nabla \cdot}
\newcommand{\Curl}{\nabla \times}
\newcommand{\mn}{{\rm min}}
\newcommand{\mx}{{\rm max}}
\newcommand{\lang}{\langle}
\newcommand{\rang}{\rangle}
\newcommand{\Md}{\partial}
\newcommand{\MP}{\bigoplus}
\newcommand{\MT}{\bigotimes}
\newcommand{\MU}{\bigcup}
\newcommand{\MI}{\bigcap}
\newcommand{\MV}{\bigvee}
\newcommand{\MW}{\bigwedge}
\newcommand{\ov}[1]{\overline{#1}}

\newcommand{\Ga}{\alpha}
\newcommand{\Gb}{\beta}
\newcommand{\Gd}{\delta}
\newcommand{\Ge}{\epsilon}
\newcommand{\Gve}{\varepsilon}
\newcommand{\Gf}{\phi}
\newcommand{\Gvf}{\varphi}
\newcommand{\Gg}{\gamma}
\newcommand{\Gc}{\chi}
\newcommand{\Gi}{\iota}
\newcommand{\Gk}{\kappa}
\newcommand{\Gl}{\lambda}
\newcommand{\Gn}{\eta}
\newcommand{\Gm}{\mu}
\newcommand{\Gv}{\nu}
\newcommand{\Gp}{\pi}
\newcommand{\Gt}{\theta}
\newcommand{\Gvt}{\vartheta}
\newcommand{\Gr}{\rho}
\newcommand{\Gvr}{\varrho}
\newcommand{\Gs}{\sigma}
\newcommand{\Gvs}{\varsigma}
\newcommand{\Gj}{\tau}
\newcommand{\Gu}{\upsilon}
\newcommand{\Go}{\omega}
\newcommand{\Gx}{\xi}
\newcommand{\Gy}{\psi}
\newcommand{\Gz}{\zeta}
\newcommand{\GD}{\Delta}
\newcommand{\GF}{\Phi}
\newcommand{\GG}{\Gamma}
\newcommand{\GL}{\Lambda}
\newcommand{\GP}{\Pi}
\newcommand{\GT}{\Theta}
\newcommand{\GS}{\Sigma}
\newcommand{\GU}{\Upsilon}
\newcommand{\GO}{\Omega}
\newcommand{\GX}{\Xi}
\newcommand{\GY}{\Psi}

\newcommand{\BGa}{\bfm\alpha}
\newcommand{\BGb}{\bfm\beta}
\newcommand{\BGd}{\bfm\delta}
\newcommand{\BGe}{\bfm\epsilon}
\newcommand{\BGve}{\bfm\varepsilon}
\newcommand{\BGf}{\bfm\phi}
\newcommand{\BGvf}{\bfm\varphi}
\newcommand{\BGg}{\bfm\gamma}
\newcommand{\BGc}{\bfm\chi}
\newcommand{\BGi}{\bfm\iota}
\newcommand{\BGk}{\bfm\kappa}
\newcommand{\BGl}{\bfm\lambda}
\newcommand{\BGn}{\bfm\eta}
\newcommand{\BGm}{\bfm\mu}
\newcommand{\BGv}{\bfm\nu}
\newcommand{\BGp}{\bfm\pi}
\newcommand{\BGt}{\bfm\theta}
\newcommand{\BGvt}{\bfm\vartheta}
\newcommand{\BGr}{\bfm\rho}
\newcommand{\BGvr}{\bfm\varrho}
\newcommand{\BGs}{\bfm\sigma}
\newcommand{\BGvs}{\bfm\varsigma}
\newcommand{\BGj}{\bfm\tau}
\newcommand{\BGu}{\bfm\upsilon}
\newcommand{\BGo}{\bfm\omega}
\newcommand{\BGx}{\bfm\xi}
\newcommand{\BMx}{\bfm\xi}
\newcommand{\BGy}{\bfm\psi}
\newcommand{\BGz}{\bfm\zeta}
\newcommand{\BGD}{\bfm\Delta}
\newcommand{\BGF}{\bfm\Phi}
\newcommand{\BGG}{\bfm\Gamma}
\newcommand{\BGL}{\bfm\Lambda}
\newcommand{\BGP}{\bfm\Pi}
\newcommand{\BGT}{\bfm\Theta}
\newcommand{\BGS}{\bfm\Sigma}
\newcommand{\BGU}{\bfm\Upsilon}
\newcommand{\BGO}{\bfm\Omega}
\newcommand{\BGX}{\bfm\Xi}
\newcommand{\BGY}{\bfm\Psi}

\newcommand{\Ca}{{\cal a}}
\newcommand{\Cb}{{\cal b}}
\newcommand{\Cc}{{\cal c}}
\newcommand{\Cd}{{\cal d}}
\newcommand{\Ce}{{\cal e}}
\newcommand{\Cf}{{\cal f}}
\newcommand{\Cg}{{\cal g}}
\newcommand{\Ch}{{\cal h}}
\newcommand{\Ci}{{\cal i}}
\newcommand{\Cj}{{\cal j}}
\newcommand{\Ck}{{\cal k}}
\newcommand{\Cl}{{\cal l}}
\newcommand{\Cm}{{\cal m}}
\newcommand{\Cn}{{\cal n}}
\newcommand{\Co}{{\cal o}}
\newcommand{\Cp}{{\cal p}}
\newcommand{\Cq}{{\cal q}}
\newcommand{\Cr}{{\cal r}}
\newcommand{\Cs}{{\cal s}}
\newcommand{\Ct}{{\cal t}}
\newcommand{\Cu}{{\cal u}}
\newcommand{\Cv}{{\cal v}}
\newcommand{\Cx}{{\cal x}}
\newcommand{\Cy}{{\cal y}}
\newcommand{\Cz}{{\cal z}}
\newcommand{\CA}{{\cal A}}
\newcommand{\CB}{{\cal B}}
\newcommand{\CC}{{\cal C}}
\newcommand{\CD}{{\cal D}}
\newcommand{\CE}{{\cal E}}
\newcommand{\CF}{{\cal F}}
\newcommand{\CG}{{\cal G}}
\newcommand{\CH}{{\cal H}}
\newcommand{\CI}{{\cal I}}
\newcommand{\CJ}{{\cal J}}
\newcommand{\CK}{{\cal K}}
\newcommand{\CL}{{\cal L}}
\newcommand{\CM}{{\cal M}}
\newcommand{\CN}{{\cal N}}
\newcommand{\CO}{{\cal O}}
\newcommand{\CP}{{\cal P}}
\newcommand{\CQ}{{\cal Q}}
\newcommand{\CR}{{\cal R}}
\newcommand{\CS}{{\cal S}}
\newcommand{\CT}{{\cal T}}
\newcommand{\CU}{{\cal U}}
\newcommand{\CV}{{\cal V}}
\newcommand{\CW}{{\cal W}}
\newcommand{\CX}{{\cal X}}
\newcommand{\CY}{{\cal Y}}
\newcommand{\CZ}{{\cal Z}}
\newcommand{\BCa}{{\bfm{\cal a}}}
\newcommand{\BCb}{{\bfm{\cal b}}}
\newcommand{\BCc}{{\bfm{\cal c}}}
\newcommand{\BCd}{{\bfm{\cal d}}}
\newcommand{\BCe}{{\bfm{\cal e}}}
\newcommand{\BCf}{{\bfm{\cal f}}}
\newcommand{\BCg}{{\bfm{\cal g}}}
\newcommand{\BCh}{{\bfm{\cal h}}}
\newcommand{\BCi}{{\bfm{\cal i}}}
\newcommand{\BCj}{{\bfm{\cal j}}}
\newcommand{\BCk}{{\bfm{\cal k}}}
\newcommand{\BCl}{{\bfm{\cal l}}}
\newcommand{\BCm}{{\bfm{\cal m}}}
\newcommand{\BCn}{{\bfm{\cal n}}}
\newcommand{\BCo}{{\bfm{\cal o}}}
\newcommand{\BCp}{{\bfm{\cal p}}}
\newcommand{\BCq}{{\bfm{\cal q}}}
\newcommand{\BCr}{{\bfm{\cal r}}}
\newcommand{\BCs}{{\bfm{\cal s}}}
\newcommand{\BCt}{{\bfm{\cal t}}}
\newcommand{\BCu}{{\bfm{\cal u}}}
\newcommand{\BCv}{{\bfm{\cal v}}}
\newcommand{\BCx}{{\bfm{\cal x}}}
\newcommand{\BCy}{{\bfm{\cal y}}}
\newcommand{\BCz}{{\bfm{\cal z}}}
\newcommand{\BCA}{{\bfm{\cal A}}}
\newcommand{\BCB}{{\bfm{\cal B}}}
\newcommand{\BCC}{{\bfm{\cal C}}}
\newcommand{\BCD}{{\bfm{\cal D}}}
\newcommand{\BCE}{{\bfm{\cal E}}}
\newcommand{\BCF}{{\bfm{\cal F}}}
\newcommand{\BCG}{{\bfm{\cal G}}}
\newcommand{\BCH}{{\bfm{\cal H}}}
\newcommand{\BCI}{{\bfm{\cal I}}}
\newcommand{\BCJ}{{\bfm{\cal J}}}
\newcommand{\BCK}{{\bfm{\cal K}}}
\newcommand{\BCL}{{\bfm{\cal L}}}
\newcommand{\BCM}{{\bfm{\cal M}}}
\newcommand{\BCN}{{\bfm{\cal N}}}
\newcommand{\BCO}{{\bfm{\cal O}}}
\newcommand{\BCP}{{\bfm{\cal P}}}
\newcommand{\BCQ}{{\bfm{\cal Q}}}
\newcommand{\BCR}{{\bfm{\cal R}}}
\newcommand{\BCS}{{\bfm{\cal S}}}
\newcommand{\BCT}{{\bfm{\cal T}}}
\newcommand{\BCU}{{\bfm{\cal U}}}
\newcommand{\BCV}{{\bfm{\cal V}}}
\newcommand{\BCW}{{\bfm{\cal W}}}
\newcommand{\BCX}{{\bfm{\cal X}}}
\newcommand{\BCY}{{\bfm{\cal Y}}}
\newcommand{\BCZ}{{\bfm{\cal Z}}}

\def\ii{{\rm i}}
\def\dd{{\rm d}}
\def\Im{{\it Im}}
\def\Re{{\it Re}}
\def\ca{{\cal A}}
\def\ct{{\cal T}}
\def\Ba{{\bf a}}
\def\Bb{{\bf b}}
\def\Bc{{\bf c}}
\def\Bd{{\bf d}}
\def\Be{{\bf e}}
\def\Bf{{\bf f}}
\def\Bg{{\bf g}}
\def\Bh{{\bf h}}
\def\Bi{{\bf i}}
\def\Bj{{\bf j}}
\def\Bk{{\bf k}}
\def\Bl{{\bf l}}
\def\Bm{{\bf m}}
\def\Bn{{\bf n}}
\def\Bo{{\bf o}}
\def\Bp{{\bf p}}
\def\Bq{{\bf q}}
\def\Br{{\bf r}}
\def\Bs{{\bf s}}
\def\Bt{{\bf t}}
\def\Bu{{\bf u}}
\def\Bv{{\bf v}}
\def\Bw{{\bf w}}
\def\Bx{{\bf x}}
\def\By{{\bf y}}
\def\Bz{{\bf z}}
\def\BA{{\bf A}}
\def\BB{{\bf B}}
\def\BC{{\bf C}}
\def\BD{{\bf D}}
\def\BE{{\bf E}}
\def\BF{{\bf F}}
\def\BG{{\bf G}}
\def\BH{{\bf H}}
\def\BI{{\bf I}}
\def\BJ{{\bf J}}
\def\BK{{\bf K}}
\def\BL{{\bf L}}
\def\BM{{\bf M}}
\def\BN{{\bf N}}
\def\BO{{\bf O}}
\def\BP{{\bf P}}
\def\BQ{{\bf Q}}
\def\BR{{\bf R}}
\def\BS{{\bf S}}
\def\BT{{\bf T}}
\def\BU{{\bf U}}
\def\BV{{\bf V}}
\def\BW{{\bf W}}
\def\BX{{\bf X}}
\def\BY{{\bf Y}}
\def\BZ{{\bf Z}}

\def\half{{\scriptstyle{1\over 2}}}
\newcommand{\eps}{\varepsilon}
\newcommand{\beq}{\begin{equation}}
\newcommand{\eeq}{\end{equation}}
\newcommand{\overliner}{\begin{eqnarray}}
\newcommand{\earr}{\end{eqnarray}}
\newcommand{\beqn}{\begin{equation*}}
\newcommand{\eeqn}{\end{equation*}}
\newcommand{\overlinern}{\begin{eqnarray*}}
\newcommand{\earrn}{\end{eqnarray*}}
\newcommand{\prt}{\partial}
\newcommand{\sg}{\sigma}
\newcommand{\fr}{\frac}
\newcommand{\kap}{\varkappa}
\newcommand{\diag}{\mbox{diag}}
\newcommand{\sign}{\mbox{sign}}
\def\l{\label}

%% file: waves_temp_interfaces_arxiv.bbl
\begin{thebibliography}{50}

\bibitem{Kittel} Kittel, C., Introduction to solid state physics, 7th Edition, 1995, John Wiley $\&$ Sons. 

\bibitem{Lekner} Lekner, J., Theory of reflection of electromagnetic and particle waves, 1987, Springer Verlag. 

\bibitem{Lurie} Lurie, K.A. An introduction to mathematical theory of dynamic materials, 2017, Springer Verlag. 

\bibitem{CD-L} Calos, C. \&  Deck-L\'eger, Z.-L., 2020. Space-time metamaterials, Part I: General concepts; Part II: Theory and applications. \textit{IEEE Transactions on Antennas and Propagation}. {\bf 68}, Issue 3, 1569--1582 (Part I), 1583--1598 (Part II).

\bibitem{Fink} Fink, M, 1993, Time-reversal mirrors. \textit{Journal of Physics D: Applied Physics}, {\bf 26}, No 9, 1333. 

\bibitem{Fink1} Bacot, V., Labousse, M., Eddi, A., Fink, M. $\&$ Fort, E., 2016. Time reversal and holography with space-time transformations, \textit{Nature Physics}, {\bf 12}, 972--977. 

\bibitem{MilMat} Milton, G.W. \& Mattei, 0., 2017. Field patterns: A new mathematical object. Proc. R. Soc. A 20160819.

\bibitem{MatMil1} Mattei, O. \& Milton, G.W. 2017. Field patterns without blow up. New J. Phys. 19 093022. 

\bibitem{MatMil2} Mattei, O. \&  Milton, G.W. 2017. Field patterns: A new type of wave with infinitely degenerate band structure. Europhys. Lett. 120(5), 54003.

\bibitem{Brun2012}
Brun, M., Jones, I.S. \& Movchan, A.B. 2012. Vortex-type elastic structured media and dynamic shielding. \textit{Proc. R. Soc. A} \textbf{468}, 3027-3046. 

\bibitem{Bertoldi}
Wang, P., Lu, L., \& Bertoldi, K. 2015. Topological phononic crystals with one-way elastic edge waves. \textit{Phys. Rev. Lett.} \textbf{115}, 104302.

\bibitem{moore}
Moore, J. E. 2010.  The birth of topological insulators, {\it Nature} {\bf 464}, 194-198.

\bibitem{pendry}
Pendry, J. B.,  Martin-Cano, D. \&  Garcia-Vidal, F. J. 2004. Mimicking surface plamons with structured
surfaces, {\it Science} {\bf 305}, 847-848.

\bibitem{hibbins}
Hibbins, A. P., Evans, B. R. \&  Sambles, J. R. 2005. Experimental verification of designer surface
plasmons, {\it Science} {\bf 308}, 670-672.

\bibitem{Zhao}  Zhao, Y., Zhou, X. \& Huang, G. 2020. Non-reciprocal Rayleigh waves in elastic gyroscopic medium. \textit{Journal of the Mechanics and Physics of Solids}, \textbf{ 143}, 104065.

\bibitem{Nassar} Nassar, H., Chen, H.,  Norris, A. N. $\&$ Huang, G. L.  2018. Quantization of band tilting in modulated phononic crystals. \textit{Phys. Rev B}, {\bf 97}, 014305 

\bibitem{LW}  Lurie, K.A.,  Weekes, S.L. 2006. Wave propagation and energy exchange in a spatio-temporal material composite with rectangular microstructure. \textit{J. Math. Anal. Appl.}, {\bf 314}, 286--310.

\bibitem{chirallattice2}
Carta, G.,  Jones, I.S., Movchan, N.V. \&  Movchan, A.B. 2019. Wave polarization and dynamic degeneracy in a chiral elastic lattice. \textit{Proc. R. Soc. A} {\bf 475}, 20190313. 

\bibitem{chirallattice1}
Carta, G., Jones, I.S., Movchan, N.V. \& Movchan, A.B. 2019. Wave characterisation in a dynamic elastic lattice: lattice flux and circulation. \textit{Phys Mesomech} {\bf 22}, 152-163.


\bibitem{GyroMS}
Nieves, M.J., Carta, G., Jones, I.S., Movchan, A.B. \& Movchan N.V. 2018. Vibrations and elastic waves in chiral multi-structures. \textit{Journal of the Mechanics and Physics of Solids}, \textbf{121}, 387-408.


\bibitem{Jones2020}
Jones, I.S., Movchan, N.V. \& Movchan, A.B. 2020. Two-dimensional waves in a chiral elastic chain: dynamic Green's matrices and localised defect modes. \textit{The Quarterly Journal of Mechanics and Applied Mathematics} \textbf{73(4)}, 305-328.

\bibitem{Engheta} Solis, D.M., Kastner, R. \& Engheta, N. 2021. Time-Varying Materials in Presence of Dispersion: Plane-Wave Propagation in a Lorentzian Medium with Temporal Discontinuity, https://arxiv.org/abs/2103.06142

\bibitem{Jones2021}
Jones, I.S., Movchan, N.V. \& Movchan, A.B. 2021. Chiral waves in structured elastic systems: dynamics of a meta-waveguide \textit{The Quarterly Journal of Mechanics and Applied Mathematics} submitted.

\bibitem{McLachlan}
McLachlan, N. W. Theory and application of Mathieu functions. 1951, Oxford University Press.








\end{thebibliography}
